\documentclass[11pt]{amsart}
\author{Daniel Johnstone}
\title[A Gelfand-Graev Formula and Stable Transfer]{A Gelfand-Graev Formula and Stable Transfer for $\text{SL}_\ell(F)$ and $\text{GL}_\ell(F)$ in the unramified case, $\ell$ an odd prime}

\usepackage{amsmath}
\usepackage{amsfonts}
\usepackage{amssymb}
\usepackage{amsthm}
\usepackage{algorithmic}
\usepackage{fancyhdr}
\usepackage{mathrsfs}

\usepackage{verbatim}

\usepackage{marvosym}

\newtheorem{thm}{Theorem}[section]
\newtheorem{prop}[thm]{Proposition}
\newtheorem{lem}[thm]{Lemma}
\newtheorem{cor}[thm]{Corollary}

\newtheorem{ass}[thm]{Assumption}

\providecommand{\bysame}{\leavevmode\hbox to3em{\hrulefill}\thinspace}

\begin{document}



\maketitle
\begin{abstract}
Let $F$ be a nonarchimedean local field of characteristic 0 with residual characteristic $p$ and let $\ell$ be an odd prime with $2\ell<p$. We establish and explicitly compute the local stable transfer factor $\Theta_\phi$ in the sense of \cite{SetT} associated to a natural $L$-embedding $\phi:{^LT}\to{^LG}$ for $G=\text{SL}_\ell$ for $\ell$ an odd prime and $T\subset G$ a maximal unramified elliptic torus defined over $F$. We also explicitly compute the associated stable transfer, answering in the affirmative the Questions A and B of \cite{SetT}. We do the same, explicitly computing the stable transfer factor $\Theta_{\widetilde{\phi}}$ and the associated stable transfer operator, in the related case of $\widetilde{\phi}:{^L\widetilde{T}}\to{^L\widetilde{G}}$ for $\widetilde{G}=\text{GL}_\ell$ and $\widetilde{T}\subset \widetilde{G}$ a maximal unramified elliptic torus defined over $F$.
\end{abstract}

\section*{Introduction}

For $G$ and $H$ reductive groups defined over a local field $F$ one can assign to a class of well-behaved $L$-homomorphisms $\phi: {^LH}\to {^LG}$ a map from $L$-packets of $H$ to $L$-packets of $G$ which may be realized as a map $\chi^H\to\chi^G=\phi(\chi^H)$ of stable (Harish-Chandra) characters. In \cite{SetT} it is proposed that this map of stable characters ought to be realized via a distribution $\Theta_\phi$ on the product $\mathfrak{A}_G(F)\times\mathfrak{A}_H(F)$ of the Steinberg-Hitchin bases of $G$ and $H$, respectively, such that the equality
\begin{equation}\label{DefiningEquation}\chi^G(a_G)=\int_{\mathfrak{A}_H(F)}\Theta_\phi(a_G,a_H)\chi^H(a_H)d\eta_H\end{equation}
holds. The distribution $\Theta_\phi$ is referred to as the stable transfer factor associated to $\phi$.

The equation (\ref{DefiningEquation}) is better understood in the following setup which uses the language of the Schwartz Kernel Theorem. Suppose that we have a space of functions $C_\phi(H)$ on $\mathfrak{A}_H(F)$ containing $C_c^\infty(\mathfrak{A}_H(F)^{\text{rss}})$ and all stable characters $\chi^H$ and a space of functions $C_\phi(G)$ on $\mathfrak{A}_H(F)$ containing $C_c^\infty(\mathfrak{A}_G(F)^{\text{rss}})$ and that we consider $\Theta_\phi\in\text{Hom}(C_\phi(G)\otimes C_\phi(H),\mathbb{C})$. Moreover, suppose that there exist linear maps $\mathfrak{L}_\phi:C_\phi(H)\to\text{Hom}(C_\phi(G),\mathbb{C})$ and $\mathfrak{R}_\phi:C_\phi(G)\to\text{Hom}(C_\phi(H),\mathbb{C})$ which satisfy
$$\left<\Theta_\phi,f\otimes g\right>=\left<\mathfrak{L}_\phi(g),f\right>=\left<\mathfrak{R}_\phi(f),g\right>$$
for all $f\in C_\phi(G)$ and $g\in C_\phi(H)$. In this setup, (\ref{DefiningEquation}) is the assertion that $\mathfrak{L}_\phi(\chi^H)=\phi(\chi^H)$ for all stable characters $\chi^H$ of $H$.

A related question posed in \cite{SetT} is whether or not, for each function $h^G\in C_c^\infty(G(F)^{\text{rss}})$, there exists a function $h^H\in C_c^\infty(H(F))$ such that
\begin{equation}\label{BigQuestions}\int_{G(F)}h^G(g)\chi^G(g)d\mu_G=\int_{H(F)}h^H(h)\chi^H(h)d\mu_H\end{equation}
for all $\chi^G=\phi(\chi^H)$. A subsequent question is whether or not there exists a corresponding function $h^H$ for all $h^G\in C_c^\infty(G(F))$ such that (\ref{BigQuestions}) holds; in \cite{SetT} the former is referred to as Question A and the latter as Question B. Though one often refers to the function $h^H$ as the stable transfer of the function $h^G$ corresponding to the $L$-homomorphism $\phi$, the function $h^H$ is not necessarily unique unless $H$ is a torus.

The analogous versions of Questions A and B on the level of Stienberg-Hitchin bases and in the language described above is, in part, whether or not there exists a linear map $\mathfrak{S}_\phi:C_\phi(G)\to C^\infty(\mathfrak{A}_H)$ such that we have $\mathfrak{R}_\phi(f)=\mathfrak{S}_\phi(f)d\eta_H$. In this language, an affirmative answer to Question A would demand that $\mathfrak{S}_\phi(f)$ is compactly supported whenever $f\in C_c^\infty(\mathfrak{A}_G(F)^{\text{rss}})$ and an affirmative answer to Question B would demand that the space $C_\phi(G)$ be large enough to contain the (normalized) stable orbital integrals of all functions $h\in C_c^\infty(G(F))$.

Analogous to the notion of endoscopic transfer of functions, stable transfer will ultimately afford a method to compare (stable) trace formulas on different groups. The notion of stable transfer is a component of Landlands' Beyond Endoscopy strategy begun in \cite{BE} which ultimately aims to establish the principle of Langlands functoriality in full generality. The stable transfer factor giving rise to a transfer of functions via the formula $f\mapsto\mathfrak{S}_\phi(f)$ for $\mathfrak{R}_\phi(f)=\mathfrak{S}_\phi(f)d\eta_H$ is the stable analogue of the Langlands-Shelstad transfer factor $\Delta$ of \cite{LS} giving rise to the endoscopic transfer of orbital integrals via
$$\text{Orb}(h)\mapsto\left(\gamma_H\mapsto \sum_\gamma\Delta(\gamma_H,\gamma)\text{Orb}(h)(\gamma)\right)$$
for an endoscopic group $H$ of $G$, a function $h\in C_c^\infty(G(F))$ and $\text{Orb}(h)$ its orbital integral, and where the above sum is taken over all elements $\gamma$ with image $\gamma_H$. Indeed, for fixed $\gamma_H\in H(F)^{\text{rss}}$ we may consider the transfer factor $\Delta(\gamma_H,\bullet)$ as a distributional kernel given by a weighted sum of point mass measures supported on the preimage of $\gamma_H$ in $G(F)$.

Few cases of a computation of either the stable transfer factor or of a formula for an associated stable transfer are known outside of the case of $G=\text{SL}_2$ and $H=T$ a maximal torus thereof. The existence of the stable transfer for all local fields of characteristic 0, with residual characteristic not equal to 2 in the nonarchimedian case, is established in \cite{SetT} giving an affirmative answer to both Questions A and B. More recently, work of \cite{YiannisGG} has, through very different methods, given a formula for the transfer in the case where $T$ is non-split. A computation of the corresponding stable transfer factor, often referred to as the Gelfand-Graev formula, can be found in \cite{GGPS}, despite the fact that several critical errors occur throughout the proof therein.

Beyond the case of $G=\text{SL}_2$, the author's thesis \cite{MyThesis} considers the case of $G=\text{SL}_\ell$ for $\ell$ an odd prime and $T\subset G$ an unramified maximal torus and computes the corresponding stable transfer factor. In this paper we give a simplified proof of the computation of the stable transfer factor in this case as well as explicitly compute the stable transfer operator, giving an affirmative answer to both Questions A and B. We also do the same for the related case of $\widetilde{G}=\text{GL}_\ell$ for $\ell$ an odd prime and $\widetilde{T}\subset \widetilde{G}$ an unramified maximal torus.

We remark that many simplifications occur in the case where $H$ is a torus which do not occur in general. The two most fundamental are that $H(F)$ is a locally compact abelian group so that harmonic analysis on $H(F)$ can be studied via classical means and that $\mathfrak{A}_H(F)$ may be identified with $H(F)$ which eliminates the need to lift functions from $\mathfrak{A}_H(F)$ to $H(F)$ in order to affirmatively answer Questions A and B. In the case $G=\text{SL}_\ell$ we even have that the torus $T(F)$ is compact, leading to an even simpler analysis, and our results in the case of $\widetilde{G}=\text{GL}_\ell$ will be seen to largely follow from the $G=\text{SL}_\ell$ case.

The above being said, despite the many aforementioned simplications owing to the fact that $H=T$ is a torus, our methods require relatively recent innovations in the explicit computation of supercuspidal character values. For our $L$-embedding ${^LT}\to {^LG}$ we denote by $\psi$ characters of the abelian group $\widehat{T}(F)$ and by $\chi_\psi$ the corresponding stable characters of $G(F)$. Though the following must be considered in the distributional setting in order for it to converge and is, in fact, not represented by a smooth function, our computation of the stable transfer factor is essentially that of the inverse Fourier transform
$$\Theta_\phi(a_G,t)=\int_{\widehat{T(F)}}\chi_\psi(a_G)\psi(t^{-1})d\nu_T.$$
For $\psi$ an admissible character, which in this case means merely that $\psi\neq \phi\circ \det$ for any $\phi\in\widehat{F^\times}$, $\chi_\psi$ is the character of a supercuspidal representation of $G(F)$; if $\widetilde{\psi}\in\widehat{\widetilde{T}(F)}$ is such that $\widetilde{\psi}|_{T(F)}=\psi$ we have that $\chi_\psi$ is the restriction to $G(F)^{\text{rss}}$ of the character of the representation $\pi_{\widetilde{\psi}}$ associated to $\widetilde{\psi}$ by \cite{Howe} upon identifying $\widehat{\widetilde{T}(F)}=\widehat{E^\times}$ for $E$ a finite unramified extension of $F$ of order $\ell$. As such, our methods require explicit character data for all such supercuspidal characters. Luckily, in our case this character data is completely known (albeit spread throughout various sources in the literature). We remark that the seemingly bizarre choice of $\text{SL}_\ell$ and $\text{GL}_\ell$ for $\ell$ a prime is related to concerns of character theory; the complexity of the character computations for $\text{GL}_n$ increases the more prime factors $n$ has.

Our main theorem for $G=\text{SL}_\ell$ is the following.
\begin{thm}\label{MAINSL}
Let $F$ be a nonarchimedean local field of characteristic 0 and residual characteristic $p$. Let $\ell$ be an odd prime and suppose $p>2\ell$. Then for $\phi:{^LT}\to{^LG}$ described in $\S$\ref{param} there exists a linear functional $\Theta_\psi\in\text{\emph{Hom}}(C_b^\infty(\mathfrak{A}_G)\otimes C_c^\infty(T(F)),\mathbb{C})$ and a linear map $\mathfrak{S}_\phi:\text{\emph{Hom}}(C_b^\infty(\mathfrak{A}_G),C_{\mu_T}^\infty(T(F)))$ for $C_b^\infty(\mathfrak{A}_G)$ and $C_{\mu_T}^\infty(T(F))$ as defined in $\S$\ref{OIsubsec} and $\S$\ref{transfersubsecSL}, respectively, satisfying
$$\mathfrak{L}_{\Theta_\phi}(\psi)=\chi_\psi d\eta_G$$
for all $\psi\in\widehat{T(F)}$ and
$$\mathfrak{R}_{\Theta_\phi}(f)=\mathfrak{S}_\phi(f)d\mu_T$$
for all $f\in C_b^\infty(\mathfrak{A}_G)$. Moreover, the space $C_b^\infty(\mathfrak{A}_G)$ contains all normalized stable orbital integrals of smooth compactly supported functions on $G(F)$ and $C_{\mu_T}^\infty(T(F))$ consists of locally constant functions supported on a compact subset of $T(F)$.
\end{thm}
Writing $\Theta_\phi$ as a distributional kernel, which is to say by fixing $a_G$ and considering $\Theta_\phi(a_G,\bullet)$ as a distribution on $T(F)$, we have 
$$\Theta_\phi(a_G,\bullet)=\left(\chi_1(a_G)+\Theta^+(a_G,\bullet)\right)d\mu_T+\delta^T(a_G)$$
for a smooth function $\Theta^+(a_G,\bullet)$ on $\mathfrak{A}_G(F)^{\text{rss}}\times T(F)^{\text{rss}}$ and $\delta^T(a_G)$ a distribution both defined in $\S$\ref{SLfamily} and $\chi_1$ the principal series character associated to the trivial representation defined in $\S$\ref{param}.

We note that since the space $C_b^\infty(\mathfrak{A}_G)$ contains all normalized stable orbital integrals of functions $h\in C_c^\infty(G(F))$ we thus have an affirmative answer to Questions A and B. We remark that we don't merely have existence of the distribution $\Theta_\phi$ and corresponding transfer operator $\mathfrak{S}_\phi$ but explicit expressions for each; see Propositions \ref{SLroundup} and \ref{SLtransferA}. Regarding the assumption $p>2\ell$, we expect our results to still hold true for $p>\ell$; we require the additional assumption in order to invoke the results of \cite{SLL} in order to make as straightforward as possible our discussion of explicit supercuspidal character values. On the other hand, we expect that the condition $p>\ell$ to be necessary as it ensures that all supercuspidal representations of all $F$-Levi subgroups of $G$ to be exhausted by the construction of \cite{Howe}.

An analogous result holds for $\widetilde{G}=\text{GL}_\ell$.
\begin{thm}\label{MAINGL}
Let $F$ be a nonarchimedean local field of characteristic 0 and residual characteristic $p$. Let $\ell$ be an odd prime and suppose $p>2\ell$. Then for $\widetilde{\phi}:{^L\widetilde{T}}\to{^L\widetilde{G}}$ described in $\S$\ref{param} there exists a linear functional $\Theta_{\widetilde{\psi}}\in\text{\emph{Hom}}(C_b^\infty(\mathfrak{A}_{\widetilde{G}})\otimes C^\infty(\widetilde{T}(F)),\mathbb{C})$ and a linear map $\widetilde{\mathfrak{S}}_{\widetilde{\phi}}:\text{\emph{Hom}}(C_b^\infty(\mathfrak{A}_{\widetilde{G}}),C_{\mu_{\widetilde{T}}}^\infty(\widetilde{T}(F)))$ for $C_b^\infty(\mathfrak{A}_{\widetilde{G}})$ and $C_{\mu_{\widetilde{T}}}^\infty(\widetilde{T}(F))$ defined in $\S$\ref{OIsubsec} and $\S$\ref{GLtransfersec}, respectively, satisfying
$$\mathfrak{L}_{\Theta_{\widetilde{\phi}}}(\widetilde{\psi})=\chi_{\widetilde{\psi}}d\eta_{\widetilde{G}}$$
for all $\psi\in\widehat{T(F)}$ and
$$\mathfrak{R}_{\Theta_\phi}(f)=\mathfrak{S}_b(f)d\mu_{\widetilde{T}}$$
for all $f\in C_b^\infty(\mathfrak{A}_{\widetilde{G}})$. Moreover, the space $C_b^\infty(\mathfrak{A}_{\widetilde{G}})$ contains all normalized stable orbital integrals of smooth compactly supported functions on $\widetilde{G}(F)$ and $C_{\mu_{\widetilde{T}}}^\infty(\widetilde{T}(F))$ consists of locally constant functions supported on a compact subset of $\widetilde{T}(F)$.
\end{thm}
This gives an affirmative answer to Questions A and B and again we have explicit descriptions of $\Theta_{\widetilde{\phi}}$ and transfer operator $\mathfrak{S}_{\widetilde{\phi}}$. Additional subtleties enter this case given that $\widetilde{T}(F)$ is not compact but we essentially proceed by relating the situation for $\widetilde{G}(F)$ to that of $G(F)$. While it is not true that the multiplication map $Z_{\widetilde{G}}(F)\times G(F)\to \widetilde{G}(F)$ is a bijection, indeed it is neither injective nor surjective, we do have that the multiplication map $Z_{\widetilde{G}}(F)\times G(F)_{0^+}\to \widetilde{G}(F)_{0^+}$ which allows us to translate many of our results from the previous case to the latter.

The paper is organized as follows.

In Section \ref{notation} we set the notation which we will need throughout the paper, notably with regards to the Fourier Transform and the Steinberg-Hitchin base, as well as normalize various measures which will appear throughout our computations.

In Section \ref{somanyprelims} we recall various notions pertaining to Bruhat-Tits buildings and the notion of depth for elements and representations of a $p$-adic group. We also use these notions to establish various facts which will facilitate computations in later sections. We also make explicit the $L$-embeddings $\phi:{^LT}\to{^LG}$ and $\widetilde{\phi}:{^L\widetilde{T}}\to{^L\widetilde{G}}$ and explicitly compute the associated stable characters $\chi_\psi$ and $\chi_{\widetilde{\psi}}$ of $G(F)$ and $\widetilde{G}(F)$, respectively.

In Section \ref{DistandOrbsec} we discuss the properties of the distributions with which we will work, notably establishing Proposition \ref{inversionprop} which is a key ingredient to our proof of Theorem \ref{MAINSL}. In short, we introduce the elementary notion of ``smooth families'' and it is via this notion we are able to glue together the individual distributions $\Theta_\phi(a_G,\bullet)$ on $T(F)$ for each $a_G$ in the regular semisimple locus of the Steinberg-Hitchin base $\mathfrak{A}_G(F)^{\text{rss}}$ into a distribution on the product $\mathfrak{A}_G(F)^{\text{rss}}\times T(F)$ and establish that this distribution indeed has the properties we desire. We also recall necessary facts about the asymptotics of (stable) normalized orbital integrals.

In Section \ref{secSL} we prove Theorem \ref{MAINSL}. We remark that Propositions \ref{SLnotT} and \ref{SLinT} essentially contain the content of the main theorem of \cite{MyThesis}.

In Section \ref{secGL} we prove Theorem \ref{MAINGL}. The proof requires some additional preliminaries which are not necessary in the case of $\text{SL}_\ell$ which largely arise from the fact that $\widetilde{T}(F)$, unlike $T(F)$, is not compact. These additional results are contained within $\S$\ref{helloGL}.

This results of this paper are a significant extension of those of \cite{SetT}. Recent breakthroughs in the explicit character computations for reductive $p$-adic groups have been a fundamental ingredient of our computations; for the scope of this paper we are able to rely primarily upon the results of \cite{Takahashi} and \cite{SLL} but formulas for more general reductive groups can be found in \cite{ASformula}. Moving beyond the scope of the present paper, we have work in progress \cite{ThesisExtension} pertaining to the case $\text{GL}_n$ and $H$ a maximal elliptic torus therein, as well as the associated $G=\text{SL}_n$ case, where $n$ is composite. To deal with this case one requires formulas found in \cite{SpiceNew} so that this case is severely complicated by the increased complexity of the supercuspidal character formulas involved. Moreover, even in the special linear case, when $n$ is composite there are infinitely many non-admissible characters of $\widetilde{T}(F)$ which contribute stable characters which are parabolically induced as opposed to supercuspidal. The character formula also becomes vastly more complicated when $\widetilde{T}$ is ramified. 

Another direction of generalization was begun in \cite{Symn} and continues in the work in progress \cite{Symn2} wherein we consider $H=\text{GL}_2$, $G=\text{GL}_{n+1}$ and $\rho=\text{Sym}^n:{^LH}\to{^LG}$. This case is vastly complicated by the fact that $H(F)$ is, in addition to being non-compact, is non-abelian. The problem was able to largely be tackled by considering the equivalence classes of maximal tori $T\subset H$ and relating the transfer from $H$ to $G$ to the well-understood transfers from each $T$ to $H$. It is our hope that the results of this paper as well as \cite{ThesisExtension} will similarly be able to be used to help establish and understand further cases of functoriality.

I would like to thank my PhD thesis advisor B\`{a}o Ch\^{a}u Ng\^{o} for first introducing me to this problem. Moreover, I would like to thank Matthew Sunohara and Patrice Moisan-Roy for the long and detailed conversations which vastly influenced the writing of this paper.

\section{Notation, Conventions and Assumptions}\label{notation}

Let $F$ be a nonarchimedean local field with $\text{char}(F)=0$ and denote by $\mathcal{O}_F$ its ring of integers, $\mathfrak{p}$ its maximal ideal, $\varpi$ a fixed choice of uniformizer of $\mathcal{O}_F$, $\mathfrak{f}=\mathcal{O}_F/\mathfrak{p}_F$ its residue field as well as $p=\text{char}(\mathfrak{f})$ and $q=|\mathfrak{f}|$. We fix an algebraic closure $\overline{F}$ of $F$ and define an additive valuation $\text{ord}:F^\times\to\mathbb{Q}$ with $\text{ord}(\varpi)=1$ and where we define $\text{ord}(0)=\infty$; we also define the absolute value $|\bullet|:\overline{F}^\times\to\mathbb{R}^\times$ via $|c|=q^{-\text{ord}(c)}$. For a finite extension of $F$ with $L\subset\overline{F}$ we denote by $\mathcal{O}_L$, $\mathfrak{p}_L$, $\mathfrak{l}=\mathcal{O}_L/\mathfrak{p}_L$ and $q_L=|\mathfrak{l}|$ its ring of integers, maximal ideal, residue field and size thereof, respectively. Moreover, we denote by $e_L$ the ramification index of $L$ over $F$ and, if $L$ is a Galois extension of $F$, $\Gamma_{L|F}$ its Galois group. 

For a locally compact abelian group $X$ equipped with Haar measure $\mu_X$ we denote by $\mathcal{F}_X$ the Fourier transform where $\mathcal{F}_X:L^1(X,\mu_X)\to C_0(\widehat{X})$ given by
$$\mathcal{F}_X(g)(\psi)=\int_Xg(x)\psi(x^{-1})$$
for all $g\in L^1(X,\mu_X)$. For a given pair $(X,\mu_X)$ we equip $\widehat{X}$ with the dual Haar measure $\nu_X$ which satisfies $\mathcal{F}_{\widehat{X}}\circ\mathcal{F}_X(g)(x)=g(x^{-1})$.

If $Y$ is any totally disconnected topological space we denote by $C_c^\infty(Y)$ the vector space of locally constant compactly supported functions on $Y$ and by $C^\infty(Y)$ the vector space of locally constant functions on $Y$.

The following facts pertaining to harmonic analysis, which we exposit as a means to fix our notation to be used throughout the paper, are well known and can be found in standard references such as \cite{Rudin}. If $X$ is a locally compact abelian group $X$ which is totally disconnected we have that $\widehat{X}$ is totally disconnected as well. In this case we have that $\mathcal{F}_X$ restricts to an isomorphism $\mathcal{F}_X:C_c^\infty(X)\to C_c^\infty(\widehat{X})$ and similarly that $\mathcal{F}_{\widehat{X}}$ restricts to an isomorphism $\mathcal{F}_X:C_c^\infty(\widehat{X})\to C_c^\infty(X)$ with ${\mathcal{F}_X}^{-1}=\mathcal{F}_{\widehat{X}}^\vee$ for $\mathcal{F}_{\widehat{X}}^\vee(\alpha)(x)=\mathcal{F}_{\widehat{X}}(\alpha)(x^{-1})$ for $\alpha\in C_c^\infty(\widehat{X})$. Moreover, for $g\in C_c^\infty(X)$ and $\alpha\in C_c^\infty(\widehat{X})$ we have the Plancherel identity
\begin{equation}\label{plancherel}\int_Xg(x)\mathcal{F}_{\widehat{X}}^\vee(\alpha)(x)d\mu_X=\int_{\widehat{X}}\mathcal{F}_X(g)(\psi)\alpha(\psi)d\nu_X\end{equation}
of which we will make heavy use throughout our main computations. Similarly, we have the Fourier Inversion Formula
\begin{equation}\label{fourierinversion}
f(x)=\int_{\widehat{X}}\int_Xf(y)\psi(xy^{-1})d\mu_Xd\nu_X
\end{equation}
valid for all $f\in C_c^\infty(X)$, $x\in X$.

At a number of points in various arguments we will have to work with quotient measures and quotient formulas for integration. If $X$ is a locally compact abelian group and $Y$ a closed subgroup we denote by $\frac{\mu_X}{\mu_Y}$ the measure on $X/Y$ which satisfies 
$$\int_Xf(x)d\mu_X=\int_{X/Y}\int_Yf(\dot{x}y)d\mu_Y\frac{d\mu_X}{d\mu_Y}$$
where, indeed, we shall write ${\displaystyle\frac{d\mu_X}{d\mu_Y}}$ instead of ${\displaystyle d\frac{\mu_X}{\mu_Y}}$. Moreover, we denote by $Y^\perp$ the subset of $\widehat{X}$ with trivial restriction to $Y$ and note that $Y^\perp$ may be naturally identified with $\widehat{X}/\widehat{Y}$ and that $\widehat{Y}$ may be naturally identified with $\widehat{X}/X^\perp$.

For $X=S(F)$ for $S$ a torus defined over $F$ we will write $\mathcal{F}_S=\mathcal{F}_{S(F)}$ and $\mu_S=\mu_{S(F)}$ as well as $\mathcal{F}_{\widehat{S}}=\mathcal{F}_{\widehat{S(F)}}$ and $\nu_S=\nu_{\widehat{S(F)}}$; indeed, this will cause no confusion in what follows since we will only be considering functions and integration on $S(F)$ as opposed to on $S(L)$ for extensions $L$ of $F$. 

For a reductive group $J$ defined over $F$ with discriminant function $D_J$ we write
$$J(F)^{\text{rss}}=\left\{j\in J(F):D_J(j)=0\right\}$$
and moreover for any subset $A\subset J(F)$ we write $A^{J,\text{rss}}=A\cap J(F)^{\text{rss}}$ or simply $A^{\text{rss}}$ when the group is understood. Notably, in an attempt to make our notations less cluttered, for an $F$-torus $S\subset J$ we write $S(F)^{\text{rss}}$ in lieu of $S(F)^{J,\text{rss}}$ as this is unlikely to cause confusion. We denote by $\mu_J$ a fixed choice of Haar measure on $J(F)$ and $\eta_J$ the measure on the $F$-points of the Steinberg-Hitchin base $\mathfrak{A}_J(F)$ as defined in \cite{FLN} and $\pi_J:J(F)\to\mathfrak{A}_J(F)$ the natural map. This measure satisfies the property that for $\pi_J:H(F)\to\mathfrak{A}_J(F)$ the natural map and $S\subset J$ a maximal torus defined over $F$ we have for $\pi_S=\pi_J|_S$ that
\begin{equation}(\pi_S)_\ast(\mu_S)=\frac{1}{W_{S,F}}|D_J|^\frac{1}{2}\left.\eta_J\right|_{\pi_S(S(F))}\end{equation}
for $W_{S,F}=N_{J(F)}(S(F))/S(F)$. Indeed, the measure on the Steinberg-Hitchin base is essentially that which appears in the Weyl Integration Formula. Moreover, we let $\pi|_S=\pi_J|_{S(F)}$.

For a positive integer $n$ we denote by $\widetilde{G}=\text{GL}_n$ and $G=\text{SL}_n$ and for any maximal torus $\widetilde{S}\subset\widetilde{G}$ we denote by $S=\widetilde{S}\cap G$ which is a maximal torus of $G$. We denote by $E\subset \overline{F}$ the unramified extension of $F$ with $[E:F]=n$ and choose a maximal torus $\widetilde{T}\subset\widetilde{G}$ with $\widetilde{T}\simeq\text{Res}_{E|F}(\mathbb{G}_m)$. For our main results, notably throughout $\S$\ref{secSL} and $\S$\ref{secGL}, we will restrict to the case where $n$ is an odd prime and write $n=\ell$. We make the following assumption on $n$ and $p$.
\begin{ass}\label{charass}
We assume $p>2n$.
\end{ass}
We remark that Assumption \ref{charass} will be used in the case $n=\ell$ an odd prime will allow us to invoke the character computations of \cite{SLL}. Moreover, Assumption \ref{charass} notably ensures that $\text{gcd}(n,p)=1$ which will be used at various points throughout our arguments below.

Regarding the Steinberg-Hitchin bases of $G$ and $\widetilde{G}$, respectively, we identify $\mathfrak{A}_G(F)\subset \mathfrak{A}_{\widetilde{G}}(F)$ as follows. For $c_1,\ldots,c_{n-1}:\widetilde{G}(F)\to F$ given by the coefficients of the characteristic polynomial
$$\text{char}(\gamma)(x)=x^n+\sum_{i=1}^{n-1}(-1)^{i-n}c_{n-i}(\gamma)x^{n-i}+(-1)^n\det(\gamma)$$
we have that the evaluation map $F[c_1,\ldots,c_{n-1}][d,d^{-1}]\to F[c_1,\ldots,c_{n-1}]$ sending $d$ to $1$ induces a map
$\mathfrak{A}_G=\text{spec}F[c_1,\ldots,c_{n-1}]\hookrightarrow \text{spec}F[c_1,\ldots,c_{n-1}][d,d^{-1}]=\mathfrak{A}_{\widetilde{G}}$ so that we may identify 
$$\mathfrak{A}_G(F)\simeq F^{n-1}\simeq \left\{1\right\}\times F^{n-1}\subset F^\times\times F^{n-1}=\mathfrak{A}_{\widetilde{G}}(F).$$
Similarly, we denote by $\det:\mathfrak{A}_{\widetilde{G}}(F)\to F^\times$ the evaluation map $\det(a_{\widetilde{G}})=\det((c_1,\ldots,c_{n-1},d))=d$. 

Since $D_{\widetilde{G}}|_{G(F)}=D_G$ and $D_{\widetilde{G}}$ is constant on (stable) conjugacy classes we define $D_G(a_G)=D_G(\gamma)$ for $a_G\in\mathfrak{A}_G(F)$ and $\gamma\in G(F)$ with $\pi_G(\gamma)=a_G$ as well as $D_{\widetilde{G}}(a_{\widetilde{G}})=D_{\widetilde{G}}(a_{\widetilde{G}})$ for $a_{\widetilde{G}}\in\mathfrak{A}_{\widetilde{G}}(F)$ and $\widetilde{\gamma}\in \widetilde{G}(F)$ with $\pi_{\widetilde{G}}(\widetilde{\gamma})=a_{\widetilde{G}}$. We define $\mathfrak{A}_G(F)^{\text{rss}}$ to be the complement of the zero locus of $D_G$ or, equivalently, as the image under $\pi_G$ of $G(F)^{\text{rss}}$. We define $\mathfrak{A}_{\widetilde{G}}(F)^{\text{rss}}$ similarly. 

The action of $Z_{\widetilde{G}}(F)$ on $\widetilde{G}(F)$ gives rise to an action of $Z_{\widetilde{G}}(F)$ on $\mathfrak{A}_{\widetilde{G}}(F)$ defined via $z\cdot(c_1,\ldots,c_{n-1},d)=(\lambda c_1,\ldots,\lambda^{n-1}c_{n-1},\lambda^n d)$ for $z=\text{diag}(\lambda,\ldots,\lambda)$, $\lambda\in F^\times$. For $z\in Z_{\widetilde{G}}(F)$ and $a_G=\pi_G(\gamma)\in \mathfrak{A}_G(F)$ we will often write $za_G=\pi_{\widetilde{G}}(z\gamma)$; this decomposition is non-unique in general but will most often be used in $\S$\ref{secGL} on a subset of $\mathfrak{A}_{\widetilde{G}}(F)$ on which we indeed have uniqueness.

In terms of measures, parallel to the decomposition of measures given by
$$\int_{\widetilde{G}}h(\widetilde{g})d\mu_{\widetilde{G}}=\int_{\widetilde{G}(F)/G(F)}\int_{G(F)}h(\dot{\widetilde{g}}g)d\mu_G\frac{d\mu_{\widetilde{G}}}{d\mu_G}$$
we may decompose the measure $\eta_{\widetilde{G}}$ via
$$\int_{\mathfrak{A}_{\widetilde{G}}(F)}f(a_{\widetilde{G}})d\eta_{\widetilde{G}}=\int_{F^\times}\int_{\mathfrak{A}_{G(F)}^d}f((c_1,\ldots,c_{n-1},d))d\eta_G^d d\mu_{F^\times}$$
for $\mathfrak{A}_{G(F)}^d=\{a_{\widetilde{G}}\in\mathfrak{A}_{\widetilde{G}}(F):\det(a_{\widetilde{G}})=d\}$ and $d\eta_G^d$ a suitable measure on $\mathfrak{A}_{G(F)}^d$. Relatedly, we normalize our measure $\mu_{Z_{\widetilde{G}}}$ on $Z_{\widetilde{G}}(F)$ so that the restriction of $\eta_{\widetilde{G}}$ to $Z_{\widetilde{G}}(F)\cdot\mathfrak{A}_G(F)$ is the pushforward of the product measure $\mu_{Z_{\widetilde{G}}}\times\eta_G$ under the action map.

\section{Character Formulas and Parameterizations}\label{somanyprelims}

In this section we will recall and make somewhat more explicit various results pertaining to the characters of supercuspidal representations of $\widetilde{G}(F)$ which arise via the construction of Howe \cite{Howe} from admissible characters $\widetilde{\psi}$ of an unramified elliptic maximal torus $\widetilde{T}(F)$ of $\widetilde{G}(F)$. 

We must work with modified notions of depth for elements of $\widetilde{G}(F)$ and characters of $\widetilde{T}(F)$ which take into account central multiplication and the Howe factorization of characters, respectively. These notions are not original; similar definitions for group elements appear in \cite{SLL}, and \cite{SL2}, for both elements and characters in \cite{DeThesis}, and notions of depth which take into account the Howe factorization of characters are implicit in the supercuspidal representation constructions of both \cite{Howe} and \cite{Construction}. This being said, we aim to both streamline and recharacterize these notions in a manner which will aid in many of our computations to follow. In addition, we perform a number of computations relating to these notions which will be used in the proofs of our main theorems.

In $\S$\ref{param} we make explicit the $L$-embeddings $^L\widetilde{T}\to{^L\widetilde{G}}$ and $^LT\to{^LG}$ with which we will work. Notably, we identify and compute the characters of the (principal series) representations of $\widetilde{G}(F)$ and $G(F)$ which correspond to the non-admissible characters of $\widetilde{T}(F)$ and $T(F)$, respectively.

\subsection{Maximal Depth for Elements}\label{maxdepthsection}

We repeat here the definitions pertaining to the notion of depth of an element of a $p$-adic group initially introduced in \cite{MP1} and \cite{MP2} though we will use the normalization of the filtration found different than which is found therein which instead may be found in, for example, \cite{ADeBMKT} and \cite{KimMurn}. With regards to this normalization, to be completely explicit, for $S=\text{Res}_{L|F}\mathbb{G}_m$ and identifying $S(F)=L^\times$, in this normalization we have for $r>0$ that $\gamma\in S(F)_r$ if and only if $\text{ord}(\gamma-1)\ge r$.

We recall that for a reductive group $J$ defined over $F$ for each $x\in\mathcal{B}(J,F)$ where $\mathcal{B}(J,F)$ is the (enlarged) Bruhat-Tits building defined in \cite{BT1} and \cite{BT2} there exists a compact open subgroup $J(F)_x$ of $J(F)$ equipped with a decreasing filtration $\{J(F)_{x,r}\}_{r\ge 0}$ (normalized as discussed above) with $J(F)_{x,0}=J(F)_x$. For $x\in\mathcal{B}(J,F)$ and writing $J(F)_{x,r^+}=\bigcup_{s>r}J(F)_s$ we define $d_{J(F),x}(\gamma)$ for each $\gamma\in J(F)$ by setting $d_{J(F),x}(\gamma)=-\infty$ if $\gamma\notin J(F)_x$ and otherwise setting $d_{J(F),x}(\gamma)=r$ if $\gamma\in J(F)_{x,r}$ but $\gamma\notin J(F)_{x,r^+}$. The depth $d_{J(F)}(\gamma)$ for $\gamma\in J(F)$ is then defined to be $-\infty$ if $\gamma\notin J(F)_x$ for any $x\in\mathcal{B}(J,F)$ and otherwise is defined to be
$$d_{J(F)}(\gamma)=\max_{x\in\mathcal{B}(J,F)}d_{J(F),x}(\gamma).$$
Moreover, the maximal depth $d_{J(F)}^+(\gamma)$ is defined to be $-\infty$ if $\gamma\notin Z_J(F)J(F)_x$ for any $x\in\mathcal{B}(J,F)$ and otherwise is defined via
$$d_{J(F)}^+(\gamma)=\max_{z\in Z_G(F)}d_{J(F)}(z\gamma).$$
We remark that, by virtue of Theorem \ref{depthprop} and Corollary \ref{slstuff}, beyond $\S$\ref{maxdepthsection} we will not need to make any subsequent reference to Bruhat-Tits buildings in any detail.

We have for a torus $S$ defined over $F$ that $S(F)_r:=S(F)_{x,r}$ for all $x\in\mathcal{B}(S,F)$ and that an explicit description of $S(F)_r$ for $r>0$ is given by
$$S(F)_r=\left\{s\in S(F):\text{ord}(\chi(s)-1)\ge r\;\forall\chi\in X^\ast(S)\right\}$$
We introduce the following definition: for a maximal torus $S\subset J$ we define $d_{S,J}(\gamma)$ for $\gamma\in S(F)$ to be $-\infty$ if $\gamma\notin Z_J(F)S(F)_0$ and for $\gamma\in Z_J(F)S(F)_0$ we define
$$d_{S,J}(\gamma)=\max\left\{r\ge 0: \text{ord}(\alpha(\gamma)-1)\ge r\;\forall\alpha\in \Phi(J,S)\right\}.$$

A number of our arguments to follow involve the relationship between depth and raising elements to powers. With regards to this, we have the following elementary fact.
\begin{lem}\label{powerlem}
For any finite extension $L$ of $F$, $c\in (L^\times)_{0^+}$ and $m$ a positive integrer, we have 
$$\text{\emph{ord}}(c^m-1)=\text{\emph{ord}}(c-1)$$
if $\text{\emph{gcd}}(p,m)=1$. Notably, Assumption \ref{charass} implies $\text{\emph{ord}}(c^m-1)=\text{\emph{ord}}(c-1)$.
\end{lem}
\begin{proof}
In this case we may write $c=1+u\varpi_L^m$ for $u\in\mathcal{O}_L^\times$ and $\varpi_L\in\mathcal{O}_L$ a uniformizer where $me_L=\text{ord}(c-1)$. We have that $c^m\in 1+mu\varpi_L^m+(\varpi_M)^{m+1}$ where $mu\in\mathcal{O}_L^\times$ since $\text{gcd}(p,m)=1$; the result follows.
\end{proof}

The first of our main facts is that the notions of the maximal depth $d^+$ and the functions $d_{S,\widetilde{G}}$ coincide in the following sense.
\begin{thm}\label{depthprop}
For $\gamma\in \widetilde{G}(F)^{\text{\emph{rss}}}$ and a maximal $F$-tamely ramified torus $\widetilde{S}\subset \widetilde{G}$ with $\gamma\in Z_{\widetilde{G}}(F)\widetilde{T}(F)_0$ we have
$$d_{\widetilde{G}(F)}^+(\gamma)=d_{\widetilde{S},\widetilde{G}}(\gamma).$$
\end{thm}
We remark that Theorem \ref{depthprop} allows us to essentially ignore any notions relating to the Bruhat-Tits building henceforth by allowing us to instead working only with valuations and roots. We expect that Theorem \ref{depthprop} holds in great generality but for our purposes we require it only for $\widetilde{G}$ (and in the form of Corollary \ref{slstuff} for $G$). To do so we require a number of lemmas which, while elementary, are somewhat involved.

\begin{lem}\label{depthsplit}
For any finite extension $L$ of $F$ and $a\in \widetilde{A}(L)_{0}$ we have that $d_{\widetilde{G}(L)}^+(a)=d_{\widetilde{A},\widetilde{G}}(a)$.
\end{lem}
\begin{proof}
By \cite{ADeBMKT} Corollary 2.2.11 we have that $d_{\widetilde{G}(L)}(za)=d_{\widetilde{A}(L)}(za)$ for all $z\in Z_{\widetilde{G}}(L)$ and hence that there exists $z_0\in Z_{\widetilde{G}}(L)$ such that $a^\prime=z_0a$ satisfies $d_{\widetilde{G}(L)}(a^\prime)=d_{\widetilde{G}(L)}^+(a^\prime)=d_{\widetilde{G}(L)}^+(a)$. It remains to compute $d_{\widetilde{G}(L)}(a^\prime)$.

We first claim that $d_{A(L)}(a^\prime)=\min_i\{\text{ord}(a^\prime_i-1)\}$. For $1\le i\le n$ let $\chi_i\in X^\ast(A)$ be given by $\chi_i(\text{diag}(c_1,\ldots,c_n))=c_i$. Since $a^\prime\in A(L)_0$ we must have $\text{ord}(a^\prime_i)=0$ for $1\le i\le n$. Moreover, if $\text{ord}(a^\prime_{i_0}-1)=0$ for some $1\le i_0\le n$ we have that $\text{ord}(\chi_{i_0}(a^\prime)-1)=\text{ord}(a^\prime_{i_0}-1)=0$ so that ${\displaystyle d_{A(L)}(a^\prime)=0=\min_{1\le i\le n}\text{ord}(a_i-1)}$, establishing the claim in this case so that we henceforth assume $\text{ord}(a_i-1)>0$ for $1\le i\le n$. For $\chi\in X^\ast(A)$ with $\chi\neq 1$ we have for $c\in A(L)$ that $\chi(c)=\prod_{i=1}^nc_i^{m_i}$ for some $m_i\in\mathbb{Z}$ for $1\le i\le n$ with $m_i\neq 0$ for at least one $1\le i\le n$. Writing $(a_i^\prime)^{m_i}=1+u_i\varpi^{k_i}$ for $u_i\in\mathcal{O}_F^\times$ for $1\le i\le n$ with $m_i\neq 0$ and $k_i=\infty$ for $m_i=0$ we have
$$\chi(a^\prime)-1=\left(\prod_{i=1}^n(a_i^\prime)^{m_i}\right)-1\in\sum_{i:m_i\neq 0}u_i\varpi^{k_i}+\varpi^{\text{min}_i\left\{k_i\right\}}\mathcal{O}_F$$
so that ${\displaystyle\text{ord}(\chi(a^\prime)-1)\ge\min_{1\le i\le n}\text{ord}((a_i^\prime)^{m_i}-1)\ge \min_{1\le i\le n}\text{ord}(a_i^\prime-1)}$. This establishes the claim.

We now claim that $d_{A,G}(a^\prime)=d_{G(L)}^+(a^\prime)$. We have $d_{A,G}(a^\prime)=\min_{i<j}\{\text{ord}(a^\prime_i-a^\prime_j)\}$ and that $d_{G(L)}^+(a^\prime)=\min_i\{\text{ord}(a^\prime_i-1)\}$ by the above claim and hence clearly have that $d_{A,G}(a^\prime)\ge d^+_{A(L)}(a^\prime)$. Letting $i_0$ be such that $d^+_{G(L)}(a^\prime)=\text{ord}(a^\prime_{i_0}-1)$, if $d_{A,G}(a^\prime)> d^+_{A(L)}(a^\prime)$ we have that $\text{ord}(a_{i_0}-a_j)>\text{ord}(a_{i_0}-1)$ for all $j\neq i_0$. It follows that for $z_{i_0}=\text{diag}(a_{i_0}^\prime,\ldots,a_{i_0}^\prime)\in Z_G(L)$ we have by the claim that 
$$d^+_{G(L)}(z_{i_0}^{-1}a^\prime)=\text{min}_i\{\text{ord}((a_{i_0}^\prime)^{-1}a_i-1)\}=\text{min}_i\{\text{ord}(a_i-a_{i_0}^\prime)\}>d^+_{A(L)}(a^\prime)$$
which contradicts the fact that $d_{A(L)}(a^\prime)=d^+_{A(L)}(a^\prime)$. It follows that $d_{A,G}(a^\prime)=d_{G(L)}^+(a^\prime)$; we are done upon noting that $d^+_{A(L)}(a)=d^+_{A(L)}(a^\prime)$.
\end{proof}

A corollary of the proof of Lemma \ref{depthsplit} is the following.
\begin{cor}\label{depthdet}
For $\gamma\in \widetilde{G}(F)^{\text{\emph{rss}}}$ and a maximal $F$-tamely ramified torus $\widetilde{S}\subset \widetilde{G}$ with $\gamma\in \widetilde{S}(F)_r$ we have that $d_{\widetilde{G}(F)}(\gamma)=\min\{\text{\emph{ord}}(\rho_1-1),\ldots,\text{\emph{ord}}(\rho_n-1)\}$ where $\rho_1,\ldots,\rho_n$ are the eigenvalues of $\gamma$ and that $\det(\gamma)\in(F^\times)_r$.
\end{cor}
\begin{proof}
Suppose $\gamma\in \widetilde{S}(F)$. There is a finite extension $L$ of $F$ over which $S$ splits; let $g\in \widetilde{G}(L)$ be such that $a=\gamma^g\in \widetilde{A}(L)$. The proof of Lemma \ref{depthsplit} shows that $d_{\widetilde{G}(L)}(\gamma)=d_{\widetilde{G}(L)}(a)=\min\{\text{ord}(\rho_1-1),\ldots,\text{ord}(\rho_n-1)\}$ and we have that $d_{\widetilde{G}(L)}(\gamma)=d_{\widetilde{G}(F)}(\gamma)$ by \cite{ADeBMKT} Corollary 2.2.4. The second statement follows easily from the first.
\end{proof}

With the above in hand we may now prove Theorem \ref{depthprop}.
\begin{proof}[Proof (of Theorem \ref{depthprop})]

Since $\widetilde{S}$ is $F$-tamely ramified there exists a tamely ramified extension $L$ of $F$ over which $\widetilde{S}$ splits and therefore an element $g\in \widetilde{G}(L)$ such that $\widetilde{S}(L)^g=\widetilde{A}(L)$ so that $\gamma^g\in \widetilde{A}(L)$. By Lemma \ref{depthsplit} we have that $d^+_{\widetilde{G}(L)}(\gamma)=d^+_{\widetilde{G}(L)}(\gamma^g)=d_{\widetilde{A},\widetilde{G}}(\gamma^g)=d_{\widetilde{S},\widetilde{G}}(\gamma)$. We clearly have $d^+_{\widetilde{G}(L)}(\gamma)\ge d^+_{\widetilde{G}(F)}(\gamma)$ so it remains to show $d^+_{\widetilde{G}(L)}(\gamma)= d^+_{\widetilde{G}(F)}(\gamma)$. If $d_{\widetilde{S},\widetilde{G}}(\gamma)=0$ we are done so henceforth we assume $d_{\widetilde{S},\widetilde{G}}(\gamma)>0$. Moreover, by replacing $\gamma$ with $z\gamma$ for $z\in Z_{\widetilde{G}}(F)$ if necessary, we may assume $d^+_{\widetilde{G}(F)}(\gamma)=d_{\widetilde{G}(F)}(\gamma)$.

Let $\rho_1,\ldots,\rho_n$ denote the eigenvalues of $\gamma$ ordered such that $\text{ord}(\rho_1-1)\le\text{ord}(\rho_j-1)$ for $2\le j\le n$. If $d_{\widetilde{S}}(\gamma)<d_{\widetilde{S},\widetilde{G}}(\gamma)$ then we have $\text{ord}(\rho_i-1)<\text{ord}\left(\frac{\rho_i}{\rho_j}-1\right)=\text{ord}(\rho_i-\rho_j)$ for all $1\le i,j\le n$. Similarly to the proof of Lemma \ref{depthsplit}, it follows that there is some $c\in\mathcal{O}_L^\times$ such that $\rho_i=c\rho_i^\prime$ with $\text{ord}(c-1)=\text{ord}(\rho_1-1)$ and $\text{ord}(\rho_i^\prime-1)>\text{ord}(\rho_1-1)$ for all $1\le i\le n$. We now have various cases.

If $\text{ord}(c-1)>0$ then $\text{ord}(c^n-1)=\text{ord}(c-1)$ by Lemma \ref{powerlem}. We have that $\det(\gamma^\prime)\in (F^\times)_0$ by Corollary \ref{depthdet} and $\det(\gamma^\prime)=c^n\prod_{i=1}^n\rho_i^\prime\in F^\times$ where $\text{ord}\left(\left(\prod_{i=1}^n\rho_i^\prime\right)-1\right)>\text{ord}(c^n-1)$. Since every element of $(F^\times)_{0^+}$ has an $n^{\text{th}}$ root in $(F^\times)_{0^+}$ by Assumption \ref{charass}, there is some $b\in F^\times$ with $\text{ord}(b-1)=\text{ord}(c-1)$ and $\text{ord}(b^{-1}c-1)>\text{ord}(c-1)$. It follows that we have $d_{\widetilde{G}(F)}(b^{-1}\gamma)>d_{\widetilde{G}(F)}(\gamma)$; this is a contradiction and completes the proof in this case.

If $\text{ord}(c-1)=0$ but there is some $b\in\mathcal{O}_F^\times$ with $\text{ord}(b^{-1}c-1)>0$ then we have $d_{\widetilde{G}(F)}(b^{-1}c)>d_T(\gamma)$ as in the previous case and we are done. Henceforth, we assume that there is no $b\in\mathcal{O}_F^\times$ with $\text{ord}(b^{-1}c-1)>0$. Again $\det(\gamma)\in (F^\times)_0$ by Corollary \ref{depthdet} so that there is some $b\in\mathcal{O}_F^\times$ with $\text{ord}(b^{-1}c^n-1)>0$. It follows that the residue field $\mathfrak{e}$ of $E$ contains the $n^{\text{th}}$ root of the image of $b$ in $\mathfrak{f}$ and hence that $L$ is not totally ramified over $F$ and that there exists an element $\tau\in\text{Aut}_F(L)$ whose image in the Galois group of $\mathfrak{l}$ over $\mathfrak{f}$ is non-trivial. For $1\le i,j\le n$, $i\neq j$ with $\rho_j=\tau(\rho_i)$ there is a root $\alpha\in\Phi(\widetilde{G},\widetilde{S})$ with $\alpha(\gamma)=\frac{\rho_i}{\rho_j}$ and we have that $\text{ord}(\alpha(\gamma)-1)=\text{ord}(\rho_i-\rho_j)=0$, contradicting that $d_{\widetilde{S},\widetilde{G}}(\gamma)>0$.
\end{proof}

With the above in hand, we further have the following analogous fact concerning $G(F)$.

\begin{cor}\label{slstuff}
For $\gamma\in G(F)^{\text{\emph{rss}}}$ there exists an $n^{\text{th}}$ root of unity $\zeta\in F^\times$ with $d_{\widetilde{G}(F)}^+(\gamma)=d_{G(F)}(\zeta^{-1}\gamma)$. Notably, $d^+_{\widetilde{G}(F)}(\gamma)=d^+_{G(F)}(\gamma)$.
\end{cor}
\begin{proof}
Let $\gamma\in G(F)^{\text{rss}}$. We have nothing to prove unless $d_{\widetilde{G}(F)}^+(\gamma)>d_{G(F)}(\gamma)$, notably that $d_{\widetilde{G}(F)}^+(\gamma)>0$. Let $z\in Z_{\widetilde{G}}(F)$ be such that $d_{\widetilde{G}(F)}(z^{-1}\gamma)=d_{\widetilde{G}(F)}^+(\gamma)$. Similarly to the proof of Theorem \ref{depthprop} in this case there exists $c\in F^\times$ with $d_{\widetilde{G}(F)}(\gamma)=\text{ord}(c-1)$ and that for $z_c=\text{diag}(c,\ldots,c)$ we have $d_{\widetilde{G}(F)}(z^{-1}\gamma)=d_{\widetilde{G}(F)}^+(\gamma)$. On the other hand, we also have
$$1=\det(\gamma)=\det(z_cz_c^{-1}\gamma)=\det(z_c)\det(z_c^{-1}\gamma)=c^n \det(z_c^{-1}\gamma).$$
If $\text{ord}(c-1)>0$ we have by Lemma \ref{powerlem} that $\text{ord}(c^n-1)=\text{ord}(c-1)$ which is a contradiction since then $\text{ord}(c^n\det(z_c^{-1}\gamma)-1)=\text{ord}(c-1)<\infty$. It follows that $\text{ord}(c-1)=0$ and that $\text{ord}(c^n-1)=\text{ord}(\det(z_c^{-1}\gamma)-1)\ge d^+(\gamma)$, the inequality following Corollary \ref{depthdet}, so that $c=\zeta b$ for some root of unity $\zeta$ and $b\in F^\times$ with $\text{ord}(b-1)\ge d_{\widetilde{G}(F)}^+(\gamma)$. It follows that for $z_\zeta\in Z_G(F)$ given by $z_\zeta=\text{diag}(\zeta,\ldots,\zeta)$ we have $d_{G(F)}(z_\zeta^{-1}\gamma)\ge d_{\widetilde{G}(F)}(z^{-1}\gamma)$ and we are done.
\end{proof}

With Theorem \ref{depthprop} and Corollary \ref{slstuff} proved we will henceforth write $d^+(\gamma)=d^+_{\widetilde{G}(F)}(\gamma)$ for any $\gamma\in \widetilde{G}(F)$ since we have shown that doing so can lead to no ambiguity. Moreover, we let $d^+(a_G)=d^+(\gamma)$ for $a_G\in\mathfrak{A}_G(F)$ and $\gamma\in G(F)$ with $\pi_G(\gamma)=a_G$ as well as $d^+(a_{\widetilde{G}})=d^+(\widetilde{\gamma})$ for $a_{\widetilde{G}}\in\mathfrak{A}_{\widetilde{G}}(F)$ and $\widetilde{\gamma}\in \widetilde{G}(F)$ with $\pi_{\widetilde{G}}(\widetilde{\gamma})=a_{\widetilde{G}}$.

We further have the following which will be of crucial importance throughout our computations in $\S$\ref{secGL}.
\begin{cor}\label{rootlem}
If $\text{\emph{gcd}}(p,n)=1$, notably if Assumption \ref{charass} holds, we have $\widetilde{G}(F)_{0^+}\subset Z_{\widetilde{G}}G(F)_{0^+}$.
\end{cor}
\begin{proof}
By Corollary \ref{depthdet} if $\gamma\in \widetilde{G}(F)_{0^+}$ we have $\det(\gamma)\in (F^\times)^{0^+}$. Since every element in $(F^\times)^{0^+}$ has an $n^{\text{th}}$ root under the assumption $\text{gcd}(p,n)=1$ it follows that there exists an element $z\in Z_{\widetilde{G}}(F)$ with $\det(z)=\det(\gamma)$ so that $\gamma=z(z^{-1}\gamma)\in Z_{\widetilde{G}}G(F)$. Moreover, we have that $z^{-1}\gamma\in Z_G(F)G(F)_{0^+}$ by Corollary \ref{slstuff} since $d^+(z^{-1}\gamma)\ge d(\gamma)>0$.
\end{proof}

Also in view of Corollary \ref{slstuff}, we define $e:G(F)\times T(F)$ via
$$e(\gamma,t)=\left\{\begin{array}{ll}
1&\text{there exists }\zeta\in Z_G(F)\text{ such that }\zeta\gamma, \zeta t\in G(F)_{0^+}\\\\
0&\text{else}
\end{array}\right.$$
where this function will appear in the definition of $\Theta_\phi$ and $\Theta_{\widetilde{\phi}}$.

We may further use the above to prove the following facts concerning the discriminant function and the sizes of various subsets of tori which will be crucial to our computations throughout $\S$\ref{secSL}. First, we have the following.
\begin{cor}\label{slelliptic}
For any maximal torus $\widetilde{S}\subset\widetilde{G}$ and $\gamma\in \widetilde{S}(F)$ we have $|D_{\widetilde{G}}(\gamma)|^{-\frac{1}{2}}\ge q^{\frac{n^2-n}{2}d^+(\gamma)}$ with equality if $\widetilde{S}=\widetilde{T}$ and $n=\ell$ an odd prime.
\end{cor}
\begin{proof}
The first statement follows by Theorem \ref{depthprop} and the definition of $d_{\widetilde{S},\widetilde{G}}$. For $\widetilde{S}=\widetilde{T}$ and $n=\ell$ we have that $\Phi(\widetilde{G},\widetilde{T})$ consists of elements of the form $\alpha_{i,j}$ for $1\le i, j\le \ell$, $i\neq j$ with, identifying $\widetilde{T}(F)=E^\times$, $\alpha_{i,j}(t)=\frac{\sigma^i(t)}{\sigma^j(t)}$ for $\sigma$ a generator of $\Gamma_{E|F}$. As such we have that $\text{ord}(\alpha_{i,j}(t)-1)=\text{ord}(\sigma^i(t)-\sigma^j(t))$ is constant for all $1\le i, j\le \ell$,  $i\neq j$. The result follows.
\end{proof}

The following technical facts, specifically c), will be integral to our proof of Theorem \ref{SLMainThm}.
\begin{lem}\label{TorusBoundLemma}
For any maximal $F$-torus $\widetilde{S}\subset\widetilde{G}$ and $S=\widetilde{S}\cap G$ we have the following.
\begin{itemize}
\item[a)] There exists a constant $Q_{\widetilde{S}}$ such that $\text{\emph{meas}}_{\mu_{\widetilde{S}}}(\widetilde{S}(F)_{m^+})= C_{\widetilde{S}}q^{-m\ell }$ for all $m\in\mathbb{Z}_{\ge 0}$.
\item[b)] There exists a constant $Q_S$ such that $\text{\emph{meas}}_{\mu_S}(S(F)_{m^+})= C_Sq^{-m(\ell-1)}$ for all $m\in\mathbb{Z}_{\ge 0}$.
\item[c)] For $Q_S$ as in b) we have, for any $F$-Levi subgroup $M\subset G$ containing $s$, that
$$\int_{S(F)_{m^+}/S(F)_{(m+1)^+}}|D_M(s)|d\mu_S\le Q_Sq^{-m\left(\frac{|\Phi_M|}{2}+\ell-1\right)}$$
where $\Phi_M$ is the absolute root system of $M$.
\end{itemize}
\end{lem}
\begin{proof}
We have that b) follows from a) since we may use the determinant map to identify $\widetilde{S}(F)/S(F)\subset F^\times$ and that there is a Haar measure $\mu_{F^\times}^{\widetilde{S}}$ on $F^\times$ satisfying
$$\int_{\widetilde{S}(F)_{m^+}}1d\mu_{\widetilde{S}}=\int_{\widetilde{S}(F)_{m^+}/S(F)_{m^+}}\int_{S(F)_{m^+}}1d\mu_S\frac{d\mu_{\widetilde{S}}}{d\mu_S}=\int_{(F^\times)_{m^+}}1d\mu_{F^\times}^{\widetilde{S}}\int_{S(F)_{m^+}}1d\mu_S$$
which establishes b) given a) since $\text{meas}_{\mu_{F^\times}^{\widetilde{S}}}((F^\times)_{m^+})$ is a constant multiple of $q^{-m}$. We now must prove a).

Any $\widetilde{S}$ is, by the existence of the rational canonical form, conjugate to a torus of the form $\widetilde{R}=\prod_{i=1}^k\widetilde{R}^i$ where $\widetilde{R}^i\simeq\text{Res}_{L^i|F}\mathbb{G}_m$ with $\sum_{i=1}^k[L^i:F]=\ell$. Further, we have $\mu_{\widetilde{R}}=\prod_{i=1}^k\mu_{\widetilde{R}^i}$ and $\widetilde{R}(F)_{m^+}=\prod_{i=1}^k\widetilde{R}^i(F)_{m^+}$. It follows that it suffices to show that there is a constant $Q_{\widetilde{R}^i}$ with $\text{meas}_{\mu_{\widetilde{R}^i}}(\widetilde{R}^i(F)_{m^+})= Q_{\widetilde{R}^i}q^{-[L^i:F] m}$ for each $i$ and, furthermore, to show each of these statements, it suffices to show that $\left|\widetilde{R}^i(F)_{m^+}/\widetilde{R}^i(F)_{(m+1)^+}\right|=q^{[L^i:F]}$. This last statement is elementary; letting $e_i$ and $f_i$ be the ramification index and inertial degree of $L^i$ over $F$, respectively, and $\mathfrak{p}_i$ be the maximal ideal of $\mathcal{O}_{L^i}$, we have $[L^i:F]=e_if_i$ and
$$\left|\widetilde{R}^i(F)_{m^+}/\widetilde{R}^i(F)_{(m+1)^+}\right|=\left|(1+\mathfrak{p}_i^{em+1})/(1+\mathfrak{p}_i^{em+e+1})\right|=(q^{f_i})^{e_i}=q^{e_if_i}=q^{[L^i:F]}.$$

For c) we first note that $|D_M(s)|^\frac{1}{2}\le q^{-m\frac{|\Phi_M|}{2}}$ for all $s\in S(F)_{m^+}$ by Corollary \ref{slelliptic}. We thus have
$$\int_{S(F)_{m^+}/S(F)_{(m+1)^+}}|D_M(s)|^\frac{1}{2}d\mu_S\le \int_{S(F)_{m^+}/S(F)_{(m+1)^+}}q^{-m\frac{|\Phi_M|}{2}}d\mu_S\le q^{-m\frac{|\Phi_M|}{2}}\text{meas}_{\mu_S}(S(F)_{m^+})$$
so that the result follows from b).
\end{proof}

\subsection{Maximal Depth for Characters}

In \cite{MP1} and \cite{MP2} the notion of the depth of an admissible representation of $J(F)$ is also introduced. When $J$ is a torus we will adopt the convention that $d(\psi)=-\infty$ when $\psi$ is the trivial representation of $J(F)$. 

For $\widetilde{\psi}\in\widetilde{T}(F)$ we define the maximal depth $d^+(\widetilde{\psi})$ via
\begin{equation}d^+(\widetilde{\psi})=d(\widetilde{\psi}|_{T(F)}).\end{equation}
We clearly have that $d^+(\widetilde{\psi})=-\infty$ whenever $\widetilde{\psi}=\phi\circ\det$ for $\phi\in \widehat{F^\times}$ but, moreover, the maximal depth has an interpretation related to the notion of the Howe factorization of $\widetilde{\psi}$ introduced in \cite{Howe}. Identifying $\widetilde{T}(F)$ with $E^\times$, therein it is shown that there exists an integer $d>0$, a tower of fields $E\supset E^0\supset\cdots\supset E^d=F$ and an increasing sequence of rational numbers $r_0<r_1<\ldots<r_{d-1}$ such that there are characters $\phi_i$ of $(E^i)^\times$ with $d(\phi_i)=r_i$ for $0\le i\le d-1$ such that
$$\widetilde{\psi}=\left(\prod_{i=0}^{d-1}\phi_i\circ N_{E|E^i}\right)\phi_d\circ N_{E|F}$$
where either $\phi_d$ is trivial or $\phi_d$ is non-trivial and $d(\phi_d)=d(\widetilde{\psi})>r_{d-1}$; see \cite{Howe} for details about additional properties satisfied by the $\phi_i$. It is easy to see that, in this notation and in the case $d^+(\widetilde{\psi})\ge 0$, we have $d^+(\widetilde{\psi})=r_{d-1}$.

The above easily implies the following.
\begin{lem}
For $n=\ell$ an odd prime we have that for all $\widetilde{\psi}\in\widehat{\widetilde{T}(F)}$ there exists $\widetilde{\psi}^\prime\in\widehat{\widetilde{T}(F)}$ and $\phi\in\widehat{F^\times}$ such that $\widetilde{\psi}=\widetilde{\psi}^\prime\omega$ for $\omega=\phi\circ\det$ and where $d(\widetilde{\psi}^\prime)=d^+(\widetilde{\psi})$.
\end{lem}

Further still, in \cite{MP1} and \cite{MP2} a filtration for the Lie algebra $\mathfrak{j}$ of a reductive group $J$: for all $x\in\mathcal{B}(J,F)$ and $r\in\mathbb{R}$ they define a decreasing filtration of lattices $\{\mathfrak{j}(F)_{x,r}\}_{r\in\mathbb{R}}$. As above, we write $\mathfrak{j}(F)_{x,r^+}=\bigcup_{s>r}\mathfrak{j}(F)_{x,s}$ and again we renormalize the filtrations as in done in \cite{ADeBMKT} and \cite{KimMurn}. Given the convention we take, for $J=\text{Res}_{E|F}\mathbb{G}_m$ for a finite extension $E$ of $F$ with ring of integers $\mathcal{O}_E$ with maximal ideal $\mathfrak{p}_E$ and ramification degree $e$ over $F$ and where we identify $\text{Lie}(J)(F)=E$, we have $X\in \text{Lie}(J)(F)_r$ if and only if $\text{ord}(X)\ge r$.

Relating the filtrations of $J(F)$ and $\mathfrak{j}(F)$ is what is often referred to as the Moy-Prasad isomorphism; if $J$ splits over a tamely ramified extension of $F$ then for any $x\in\mathcal{B}(J,F)$ and $0<s\le r\le 2s$ we have by \cite{Construction} that
$$J(F)_{x,s}/J(F)_{x,r}\simeq \mathfrak{j}(F)_{x,s}/\mathfrak{j}(F)_{x,s}$$
Moreover, in \cite{MP2} it is shown that for all $x\in\mathcal{B}(J,F)$ and all $r\in\mathbb{R}$ that 
$$\mathfrak{j}(F)_{x,r}/\mathfrak{j}(F)_{x,r^+}\simeq\mathfrak{j}(F)_{x,0}/\mathfrak{j}(F)_{x,0^+}\simeq \text{Lie}(\mathbf{J}_x)(\mathfrak{f})$$
where $\mathbf{J}_x$ is a reductive group defined over $\mathfrak{f}$ where we have $\mathbf{J}_x(\mathfrak{f})=J(\mathfrak{f})$ if $J$ is defined over $\mathcal{O}_F$.

To perform our calculations to follow, notably throughout $\S$\ref{secSL}, we require explicit computations regarding the size of various subsets of $\widehat{T(F)}$ depending on depth. We have the following.

\begin{lem}\label{CharCount}
We have for all positive integers $d$ that 
$$\left|\{\psi\in\widehat{T(F)}:d(\psi)\le 0\}\right|=|T(\mathfrak{f})|q^d=\frac{q^n-1}{q-1}q^{dn}$$
and that $|\{\psi\in\widehat{T(F)}:d(\psi)=0\}|=|T(\mathfrak{f})|-1=\frac{q^n-q}{q-1}$ as well as
$$|\{\psi\in\widehat{T(F)}:d(\psi)= d\}|=|T(\mathfrak{f})|(q^n-1)q^{(d-1)(n-1)}=\frac{q^n-1}{q-1}(q^{n-1}-1)q^{(d-1)(n-1)}$$
\end{lem}
\begin{proof}
Since $T$ is unramified it follows from Corollary \ref{depthdet} that $T(F)_{d^+}=T(F)_{d+1}$ for all integers $d\ge 0$. The definition of depth gives us that $\{\psi\in\widehat{T(F)}:d(\psi)\le d\}\simeq T(F)/T(F)_{d^+}$ and we have
$$|T(F)/T(F)_{d^+}|=|T(F)/T(F)_1|\prod_{i=1}^d|T(F)_i/T(F)_{i+1}|$$
We have that $|T(F)/T(F)_1|=|T(\mathfrak{f})|=\frac{|\mathfrak{e}^\times|}{\mathfrak{f}^\times}|=\frac{q^n-1}{q-1}$. Moreover, the Moy-Prasad isomoprhism gives us that 
$$|T(F)_i/T(F)_{i+1}|=|\mathfrak{t}(F)_i/\mathfrak{t}(F)_{i+1}|=|\mathfrak{t}(\mathfrak{f})|=q^{n-1}$$
and hence that $|T(F)/T(F)_{d^+}|=(q+1)q^{d(n-1)}$. The other claims follow.
\end{proof}

\subsection{Parametrization of Stable Characters}\label{param}

Here we make explicit the $L$-embeddings $\phi:{^LT}\to{^LG}$ and $\widetilde{\phi}:{^L\widetilde{T}}\to{^L\widetilde{G}}$ with which we will work. Throughout $\S$\ref{param} we restrict to the case $n=\ell$ an odd prime.

Since $\widetilde{T}\simeq\text{Res}_{E|F}(\mathbb{G}_m)$ we we may identify $^L\widetilde{T}\simeq\text{Ind}_{W_E}^{W_F}\mathbb{C}\rtimes\Gamma_{E|F}$ and that since $\widetilde{G}$ is split over $F$ we may identify $^L\widetilde{G}=\widetilde{G}(\mathbb{C})\times\Gamma_{E|F}$. We define $\widetilde{\phi}:{^L\widetilde{T}}\to{^L\widetilde{G}}$ via
$$(x_1,\ldots,x_\ell)\rtimes 1\mapsto\begin{pmatrix}x_1&&\\&\ddots&\\&&x_\ell\end{pmatrix}, 1\rtimes\sigma\mapsto\begin{pmatrix}&I_{\ell-1}\\1&\end{pmatrix}$$
where $I_{\ell-1}$ is the $(\ell-1)\times(\ell-1)$ identity matrix. By \cite{LocLang} we may identify the set of Langlands parameters of ${^L{\widetilde{T}}}$ with $\widehat{\widetilde{T}(F)}$ in such a way that we have $\phi\circ\widetilde{\psi}=\text{Ind}_{W_E}^{W_F}\widetilde{\psi}$. We have the following.
\begin{prop}\label{parameterization}
If $n=\ell$ is an odd prime and $\widetilde{\psi}\in\widehat{\widetilde{T}(F)}$ is admissible in the sense of \cite{Howe} we have that the parameter $\phi\circ\widetilde{\psi}$ of ${^L\widetilde{G}}$ corresponds to the representation $\pi_{\widetilde{\psi}}$ as defined in \cite{Howe}. Moreover, for any $n$, $\widetilde{\psi}$ is the central character of $\pi_{\widetilde{\psi}}$ and $\pi_{\widetilde{\psi}(\phi\circ\det)}=(\phi\circ\det)\otimes\pi_{\widetilde{\psi}}$ for all $\phi\in\widehat{F^\times}$.
\end{prop}
\begin{proof}
The first statement is \cite{SLL} Theorem 9.2. The other statements are clear from the construction of \cite{Howe}.
\end{proof}

If $\widetilde{\psi}=\omega\circ N_{E|F}$ we have by direct computation that
$$\text{Ind}_{W_E}^{W_F}\omega\circ N_{E|F}=(\omega\circ N_{E|F})\otimes \text{Ind}_{W_E}^{W_F}1=(\omega\circ N_{E|F})\otimes(1\oplus\text{sgn}_{E|F}\oplus\cdots\oplus \text{sgn}_{E|F}^{\ell-1})$$
for $\text{sgn}_{E|F}$ the character of $F^\times$ trivial on $N_{E|F}(E^\times)$ attached to $E$ via class field theory. Indeed, for $\text{Ind}_{W_E}^{W_F}1$ we have bases $e_i=\sigma^i\mathbf{1}_{W_E}$ and $v_j=\sum_{i=0}^{\ell-1}\zeta^{ij}e_i$ for $\zeta=\text{sgn}(\sigma)$ a non-trivial $\ell^{\text{th}}$-root of unity of we observe that $\sigma\cdot v_j=\text{sgn}(\sigma)^jv_j$. It is well known, see for example (Loc lang ref), that $\text{Ind}_{W_E}^{W_F}(\omega\circ N_{E|F})$ corresponds to the representation $(\omega\circ\det)\otimes\text{Ind}_{\widetilde{B}(F)}^{\widetilde{G}(F)} 1\oplus\text{sgn}_{E|F}\oplus\cdots\oplus \text{sgn}_{E|F}^{\ell-1}$ which by \cite{vanDijk} has character supported on $\widetilde{A}(F)^{\widetilde{G}(F)}$. For $\gamma\in\widetilde{A}(F)$ with $\gamma=\text{diag}(\gamma_1,\ldots,\gamma_\ell)$ we have
$$\chi_{\omega\circ N_{E|F}}(\gamma)=((\omega{\text{sgn}_{E|F}}^{-1})\circ\det)(\gamma)\frac{1}{|D_G(\gamma)|^\frac{1}{2}}u_1(\gamma)$$
for
$$u_1(\gamma)=\sum_{\tau\in\mathfrak{S}_\ell}\text{sgn}_{E|F}(\gamma_{\tau(1)})\cdots\text{sgn}_{E|F}^\ell(\gamma_{\tau(\ell)})$$
The support of the function $u_1$ will become an important ingredient in our computations in $\S$\ref{secGL}. We have the following.
\begin{lem}\label{u1lem}
For $\gamma\in \widetilde{G}(F)^{\text{\emph{rss}}}$ we have $u_1(\gamma)=0$ unless $\det(\gamma)\in N_{E|F}(E^\times)$.
\end{lem}
\begin{proof}
Let $\gamma=\text{diag}(\gamma_1,\ldots,\gamma_\ell)$ and write $\text{sgn}_{E|F}(\gamma_i)=\zeta^{m_i}$ for $\zeta$ an $\ell^{\text{th}}$ root of unity. Also, writing $\gamma_{\ell+1}=\gamma_1$ and $m_{\ell+1}=m_\ell$ for notational convenience, let $\gamma^\prime=\text{diag}(\gamma_{2},\ldots,\gamma_{\ell+1})$. We clearly have $u_1(\gamma)=u_1(\gamma^\prime)$ by the definition of $u_1$. On the other hand, direct computation shows that $\text{sgn}_{E|F}(\det(\gamma))u_1(\gamma^\prime)=u_1(\gamma)$ since we have
$$\text{sgn}_{E|F}(\det(\gamma))u_1(\gamma^\prime)=\sum_{\tau\in\mathfrak{S}_\ell}\prod_{j=1}^\ell\zeta^{(j+1)m_{\tau(j+1)}}=\sum_{\tau\in\mathfrak{S}_\ell}\prod_{j=1}^\ell\zeta^{jm_{\tau(j)}}=u_1(\gamma).$$
\end{proof}
Lemma \ref{u1lem} implies that we have
$$\chi_{\omega\circ N_{E|F}}(\gamma)=(\omega\circ\det)(\gamma)\frac{1}{|D_G(\gamma)|^\frac{1}{2}}u_1(\gamma)$$
since $u_1(\gamma)\neq 0$ implies $\det(\gamma)\in N_{E|F}(E^\times)$.

The map $\widetilde{\phi}:{^L\widetilde{T}}\to{^L\widetilde{G}}$ induces a map $\phi:{^LT}\to{^LG}$ by identifying $^LT={^L\widetilde{T}}/Z_{\widetilde{G}}(\mathbb{C})$ and $^LG={^L\widetilde{G}}/Z_{\widetilde{G}}(\mathbb{C})$. By \cite{SLLpack} Theorem 4.1 and again by \cite{SLL} Theorem 9.2 we have that the stable characters of $G(F)$ are the restrictions to $G(F)$ of characters of representations of $\widetilde{G}(F)$ and that for $\psi\in \widehat{T(F)}$ we have
$$\chi_\psi=\left.\chi_{\widetilde{\psi}}\right|_{G(F)^{\text{rss}}}$$
for any $\widetilde{\psi}\in\widehat{\widetilde{T}(F)}$ with $\left.\widetilde{\psi}\right|_{T(F)}=\psi$.

\subsection{The Local Character Expansion}\label{LCEsubsection}

Here we write down a very explicit formula for the Local Character Expansion for the representations $\pi_{\widetilde{\psi}}$. While other results are also necessary, the main computations we rely upon are those found in \cite{LCE}. The following is well known.
\begin{prop}
For $\widetilde{\psi}\in\widehat{\widetilde{T}(F)}$ which is admissible in the sense of \cite{Howe} and $\gamma\in \widetilde{G}(F)$ with $d^+(\gamma)>d^+(\widetilde{\psi})$ we have, writing $\gamma=z\gamma^\prime$ for $z\in Z_{\widetilde{G}}(F)$ and $\gamma^\prime\in G(F)$, that
$$\chi_{\widetilde{\psi}}(\gamma)=\text{\emph{LCE}}_{\widetilde{\psi}}(\gamma)=\widetilde{\psi}(z)\sum_{\mathcal{O}\le\mathcal{O}_\gamma}c_\mathcal{O}(\pi_{\widetilde{\psi}}) \widehat{\mu}_{\mathcal{O}}(\gamma^\prime-1)$$
\end{prop}
\begin{proof}
It is known by \cite{DeBLCE} Theorem 3.5.2, that Assumption \ref{charass} allows us to replace $\log(\gamma)$ with $\gamma-1$ follows from \cite{DeThesis} $\S$3.7, and the last statement of Proposition \ref{parameterization} allows us to work with $d^+(\widetilde{\psi})$ instead of $d(\widetilde{\psi})=d(\pi_{\widetilde{\psi}})$.
\end{proof}

Our present goal is to make more explicit the formula for $\text{LCE}_{\widetilde{\psi}}(\gamma)$ and as such require formulas for both the $c_\mathcal{O}(\pi_{\widetilde{\psi}})$ and the $\widehat{\mu}_{\mathcal{O}}$. The following formula, the notation within which will be defined below, holds; the remainder of $\S$\ref{LCEsubsection} will be devoted to the proof and explication of Proposition \ref{LCEexplicit}.
\begin{prop}\label{LCEexplicit}
For $\gamma=z\gamma^\prime$ we have that
$$\text{\emph{LCE}}_{\widetilde{\psi}}(\gamma)=\widetilde{\psi}(z)\sum_{\mathcal{O}\le\mathcal{O}_\gamma}C_\mathcal{O}q^{\frac{|\Phi_\mathcal{O}|}{2}d^+(\widetilde{\psi})}\sum_{w\in W_\mathcal{O}(\gamma)}\frac{|D_{M_\mathcal{O}}(\gamma^w)|^\frac{1}{2}}{|D_{\widetilde{G}}(\gamma)|^\frac{1}{2}}$$
for
$$C_\mathcal{O}=\frac{n(-1)^{n+r_\mathcal{O}}(r_\mathcal{O}-1)!}{w_\mathcal{O}}.$$
\end{prop}

We require some facts about the set $\mathcal{N}$ of nilpotent orbits in $\mathfrak{gl}_n(F)$ which we obtain from \cite{Nilpotent}. To each standard $F$-Levi subgroup $M\subset \widetilde{G}$ one may associate a nilpotent orbit $\mathcal{O}_M\in\mathcal{N}$ and we have that for all $\mathcal{O}\in\mathcal{N}$ there is some standard $F$-Levi subgroup $M$ for which $\mathcal{O}=\mathcal{O}_M$. As we will notationally require it, for any $\mathcal{O}\in\mathcal{N}$ we denote by $M_\mathcal{O}$ the standard $F$-Levi subgroup for which $\mathcal{O}=\mathcal{O}_{M_\mathcal{O}}$. Also, for any standard Levi subgroup $M$ of $\widetilde{G}$ there exist integers $r_M\ge 1$ and $m_{M,i}$ for $1\le i\le r_M$ such that $M=\prod_{i=1}^{r_M}\text{GL}_{m_{M,i}}$ and $\sum_{i=1}^{r_M}m_{M,i}=n$; we write $\vec{m}_M=(m_{M,1},\ldots, m_{M,r_M})$ so that $\vec{m}_M$ is a partition of $n$ which corresponds to $M$. Moreover, if $M=M_\mathcal{O}$ we write $r_\mathcal{O}=r_{M_\mathcal{O}}$ and denote by $\Phi_\mathcal{O}$ the root system of $M_\mathcal{O}$.

For $\gamma\in \widetilde{G}(F)^{\text{rss}}$ we define $\mathcal{O}_\gamma=\mathcal{O}_{M_\gamma}$ for $M_\gamma$ the centralizer of the split component of $C_{\widetilde{G}}(\gamma)$. For standard Levi subgroups $M,M^\prime$ of $\widetilde{G}$ with split components $A_M$ and $A_{M^\prime}$, respectively, the set $W(A_M,A_{M^\prime})$ is defined in \cite{LCE} Lemma 5.1 to be 
$$W(A_M,A_{M^\prime})=\{s:A_M\to A_{M^\prime}:s\text{ is an injection, } s(a)=yay^{-1}\text{ for some y}\in\widetilde{G}(F)\}$$
and we let $W_{\mathcal{O}}(\gamma)=W(A_{M_\mathcal{O}},A_{M_\gamma})$. In this notation, \cite{LCE} Lemma 5.1 gives that
$$\widehat{\mu}_{\mathcal{O}}(\gamma-1)=\frac{1}{|D_G(\gamma)|^\frac{1}{2}}\sum_{s\in W_{\mathcal{O}}(\gamma)}|D_M({^s\gamma})|^\frac{1}{2}$$
for $\gamma\in\widetilde{G}(F)_{0^+}$. Moreover, denoting by $w_\mathcal{O}$ the size of the group $N_{\widetilde{G}}(A_{M_\mathcal{O}})/M_\mathcal{O}$, \cite{LCE} Proposition 5.3 gives that
$$c_\mathcal{O}(\pi_{\widetilde{\psi}})=(-1)^{n+r_\mathcal{O}}n(r_\mathcal{O}-1)!\frac{1}{w_\mathcal{O}}q^{\frac{n^2-n-\text{dim}(\mathcal{O})}{2}d^+(\widetilde{\psi})}$$
so that we will have proven Proposition \ref{LCEexplicit} upon establishing the following explicit computation of $\text{dim}(\mathcal{O})$.
\begin{lem}
For any $\mathcal{O}\in\mathcal{N}$ we have
$$\text{\emph{dim}}(\mathcal{O})=n^2-n-|\Phi_{\mathcal{O}}|.$$
\end{lem}
This fact is well-known but as we could not find an explicit reference stating the result in this precise form we include a short proof.
\begin{proof}
For a partition $\vec{m}=(m_1,\ldots,m_r)$ of $n$ we have by \cite{Nilpotent} Corollary 6.1.4 that
$$\text{dim}(\mathcal{O}_{\vec{m}})=n^2-\sum_{i=1}^{r^t}m_i^t$$
where $\vec{m}^t=(m_1^t,\ldots,m_r^t)$ is the partition of $n$ dual to $\vec{m}$. Moreover, by \cite{Nilpotent} Theorem 7.2.3 we have for any standard $F$-Levi subgroup $M$ that $\mathcal{O}_M=\mathcal{O}_{m_M^t}$. Since $\vec{m}=(\vec{m}^t)^t$ and we easily see that
$$|\Phi_\mathcal{O}|=-n+\sum_{i=1}^{r_M}m_{M,i}^2$$
we obtain
\begin{align*}
\text{dim}(\mathcal{O}_M)=&n^2-\sum_{i=1}^{r_M}m_{M,i}^2=n^2-n-\left(-n+\sum_{i=1}^{r_M}m_{M,i}^2\right)=n^2-n-|\Phi_M|
\end{align*}
as required.
\end{proof}

\subsection{Supercuspidal Character Formulas}\label{charformsubsec}

With the above in hand we now have the following explicit character tables for the $\pi_{\widetilde{\psi}}$ for admissible $\widetilde{\psi}$. These results are known but appear in various disparate sources and for various reasons are not always written down as explicitly as we require; we include Proposition \ref{CharacterTable} to reference throughout our computations in subsequent sections.
\begin{prop}\label{CharacterTable}
Let $\gamma\in \widetilde{G}(F)^{\text{\emph{rss}}}$ and $\widetilde{\psi}\in\widehat{\widetilde{T}(F)}$ be admissible. For $\gamma\in\widetilde{T}(F)$ we have
$$\Theta_{\widetilde{\psi}}(\gamma)=q^{\frac{\ell^2-\ell}{2}\text{\emph{min}}\{d^+(\gamma), d^+(\psi)\}}\sum_{\sigma\in\Gamma}\psi(\gamma^\sigma).$$
If $\gamma\notin\widetilde{T}(F)^{\widetilde{G}(F)}$ and $d^+(\gamma)>d^+(\psi)$ we may write $\gamma=z\gamma^\prime$ for $z\in Z_{\widetilde{G}}(F)$ and $\gamma^\prime\in G(F)$ and we have
$$\Theta_{\widetilde{\psi}}(z\gamma^\prime)=\text{\emph{LCE}}_{\widetilde{\psi}}(\gamma)$$
for $\text{\emph{LCE}}_{\widetilde{\psi}}(\gamma)$ as defined in Proposition \ref{LCEexplicit}. and if $\gamma\notin\widetilde{T}(F)^{\widetilde{G}(F)}$ and $d^+(\gamma)\le d^+(\psi)$ we have $\Theta_{\widetilde{\psi}}(\gamma)=0$.
\end{prop}
\begin{proof}
The cases where $\gamma$ elliptic is elliptic follow from \cite{Takahashi} Theorem 2.2.3. For $\gamma$ non-elliptic with $d^+(\gamma)>d^+(\psi)$ the result is that of Proposition \ref{LCEexplicit}; we note that if $d^+(\gamma)>d^+(\psi)$ then the formulas of \cite{Takahashi} and those of Proposition \ref{LCEexplicit} agree. For non-elliptic $\gamma$ with $d^+(\gamma)\le d^+(\psi)$ the result follows from \cite{SLL} Corollary 6.4. 
\end{proof}

We remark that the character formula of \cite{ASformula} can be used to easily show that $\Theta_{\widetilde{\psi}}(\gamma)=0$ for nonelliptic $\gamma$ with $d^+(\gamma)< d^+(\psi)$ but cannot easily handle the case when $d^+(\gamma)= d^+(\psi)$ which is why we have instead invoked \cite{SLL}. We further remark that the formulas of \cite{Takahashi} and \cite{SLL} agree for $\gamma\in\widetilde{T}(F)$, $d^+(\gamma)\le d^+(\psi)$.

Henceforth we define $\chi_\psi(a_G)=\chi_\psi(\gamma)$ for $a_G\in\mathfrak{A}_G(F)$ and $\gamma\in G(F)$ with $\pi_G(\gamma)=a_G$ for all $\psi\in\widehat{T(F)}$. We define $\chi_{\widetilde{\psi}}(a_{\widetilde{G}})$ similarly for $a_{\widetilde{G}}\in\mathfrak{A}_{\widetilde{G}}(F)$, $\widetilde{\psi}\in\widehat{\widetilde{T}(F)}$. Moreover, we remark that for an $F$-Levi subgroup $\widetilde{M}\subset\widetilde{G}$ with discriminant function $D_{\widetilde{M}}$ the map $\gamma\mapsto D_{\widetilde{M}}(\gamma)$ is not stable in general. To remedy this, we define
$$B^\mathcal{O}(\gamma)=\sum_{w\in W_\mathcal{O}(\gamma)}\frac{|D_{\widetilde{M}}(\gamma^w)|^\frac{1}{2}}{|D_{\widetilde{G}}(\gamma)|^\frac{1}{2}}$$
which is stable so that we may further define $B^\mathcal{O}(a_{\widetilde{G}})=B^\mathcal{O}(\gamma)$ for any $\gamma\in\widetilde{G}(F)$ with $\pi_{\widetilde{G}}(\gamma)=a_{\widetilde{G}}$.

\section{Distributions and Orbital Integrals}\label{DistandOrbsec}

To perform our computations to follow we must define and analyze a number of different spaces of distributions. The stable transfer factor we seek to understand is, in addition to being defined as a distribution, build up in some sense from various inverse Fourier Transforms which exist as distributions but which do not converge in the classical sense. Moreover, to properly analyze Langlands' Question B, we will also need to recall some facts pertaining to the asymptotics of orbital integrals viewed as functions on the Steinberg-Hitchin base.

\subsection{Distributions and a Schwartz Kernel Theorem}\label{distsubsec}

The facts we recall pertaining to functional and harmonic analysis as well as to the theory of distributions are well known and are recalled herein largely in order to set notation to be used henceforth; see \cite{Rudin} and \cite{Distr}.

For a complex vector space $V $ we denote by $V^\#$ its linear dual space $\text{Hom}(V,\mathbb{C})$. For $X$ a locally compact totally disconnected topological space we set $\mathcal{D}(X)=C_c^\infty(X)^\#$ and for $f\in C_c^\infty(X)$ and $\lambda\in\mathcal{D}(X)$ we write $\left<\lambda,f\right>=\lambda(f)$. This bracket notation will clarify many of our manipulations to follow, notably those involving the Fourier Transform.

For a Radon measure $\mu_X$ on $X$ and $h\in C^\infty(X)$ we define $hd\mu_X\in\mathcal{D}(X)$ via
$$\left<hd\mu_X,f\right>=\int_Xh(x)f(x)d\mu_X$$
for all $f\in C_c^\infty(X)$. Indeed, the product $x\mapsto h(x)f(x)$ is bounded and compactly supported and hence is integrable. 

If the space $X$ is, moreover, a locally compact abelian group, for any $\lambda\in\mathcal{D}(\widehat{X})$ we may define $\mathcal{F}_X^\vee(\lambda)\in\mathcal{D}(X)$ via
$$\left<\mathcal{F}_X^\vee(\lambda),f\right>=\left<\lambda,\mathcal{F}_X(f)\right>$$
for all $f\in C_c^\infty(X)$. The following lemma is elementary to prove but we include it so that it may be explicitly referenced in our computations.
\begin{lem} The following results hold.
\begin{itemize}
\item[a)] For $\alpha\in C_c^\infty(\widehat{X})$ we have that $\mathcal{F}_X^\vee(\alpha d\nu_X)=\mathcal{F}_X^\vee(\alpha)d\mu_X$.
\item[b)] For $x\in X$ and $\beta_x\in C^\infty(\widehat{X})$ given by $\beta_x(\psi)=\psi(x)$ we have that $\mathcal{F}_X^\vee(\beta_x)=\delta_x$ for $\left<\delta_x,f\right>=f(x)$ for all $f\in C_c^\infty(X)$.
\end{itemize}
\end{lem}
\begin{proof}
The first statement follows immediately from the Plancherel identity (\ref{plancherel}) and the second follows from Fourier inversion.
\end{proof}

Since we are not endowing our function spaces with any topology we have that the Tensor-Hom adjunction $(V\otimes W)^\#\simeq\text{Hom}(V,W^\#)$ yields the following simple version of the Schwartz Kernel Theorem which is all we will require to perform our computations.
\begin{prop}\label{SK}
For any locally compact totally disconnected topological spaces $X$ and $Y$ and $\lambda\in\mathcal{D}(X\times Y)$ there exist unique linear maps $\mathfrak{L}_\lambda:C_c^\infty(X)\to \mathcal{D}(Y)$ and $\mathfrak{R}_\lambda:C_c^\infty(Y)\to \mathcal{D}(X)$ satisfying
$$\left<\lambda,f\otimes g\right>=\left<\mathfrak{L}_\lambda(g),f\right>=\left<\mathfrak{R}_\lambda(f),g\right>$$
for all $f\in C_c^\infty(X)$ and $g\in C_c^\infty(Y)$. Conversely, for any map $\mathfrak{L}:C_c^\infty(X)\to \mathcal{D}(Y)$ there exists a unique $\lambda_{\mathfrak{L}}\in\mathcal{D}(X\times Y)$ satisfying 
$$\left<\lambda_{\mathfrak{L}},f\otimes g\right>=\left<\mathfrak{L}(f),g\right>$$
and for any map $\mathfrak{R}:C_c^\infty(Y)\to \mathcal{D}(X)$ there exists a unique $\lambda_{\mathfrak{R}}\in\mathcal{D}(X\times Y)$ satisfying 
$$\left<\lambda_{\mathfrak{R}},f\otimes g\right>=\left<\mathfrak{R}(g),f\right>.$$
\end{prop}

For locally compact totally disconnected topological spaces $X$ and $Y$ we call a map $\mathcal{P}:X\to\mathcal{D}(Y)$ a smooth family if for each $g\in C_c^\infty(Y)$ we have for $\mathcal{P}_g(x)=\left<\mathcal{P}(x),g\right>$ that $\mathcal{P}_g\in C_c^\infty(X)$. In the sense of part a) of the Proposition \ref{inversionprop} below, we may assign to any smooth family $\mathcal{P}$ a distribution on $X\times Y$. Proposition \ref{inversionprop} is the technical result we require to prove the main results of the paper.
\begin{prop}\label{inversionprop}
Let $X$ and $Y$ be locally compact totally disconnected topological spaces.
\begin{itemize}
\item[a)] If $\mathcal{P}:X\to\mathcal{D}(Y)$ is a smooth family there exists a unique $\lambda_\mathcal{P}\in\mathcal{D}(X\times Y)$ satisfying
$$\left<\mathfrak{L}_{\lambda_\mathcal{P}}(g),f\right>=\int_X\mathcal{P}_g(x)f(x)d\mu_X.$$
\item[b)] Suppose $Y$ is a locally compact abelian group and $h\in C^\infty(X\times\widehat{Y})$ where we write $h_x(\psi)=h(x,\psi)$ for each $x\in X$, $\psi\in\widehat{Y}$. The map $\mathcal{P}^h:X\to\mathcal{D}$ given by $\mathcal{P}^h(x)=\mathcal{F}_Y^\vee(h_x)$ is a smooth family.
\item[c)] Maintaining the notation of b) and writing $\lambda=\lambda_{\mathcal{P}^h}$, if $Y$ is compact we have that
$$\mathfrak{L}_\lambda(\psi)=h(\bullet,\psi)d\mu_X$$
for all $\psi\in\widehat{Y}$.
\end{itemize}
\end{prop}
\begin{proof}
By Proposition \ref{SK} the formula in a) defines a unique element of $\mathcal{D}(X\times Y)$. Indeed, existence is clear, and uniqueness follows since if $\lambda$ and $\lambda^\prime$ are such that
$$\int_X\mathcal{P}_g(x)f(x)d\mu_X=\left<\mathfrak{L}_{\lambda}(g),f\right>=\left<\mathfrak{L}_{\lambda^\prime}(g),f\right>$$
for all $f\in C_c^\infty(X)$ and $g\in C_c^\infty(Y)$ we must have that $\mathfrak{L}_{\lambda}(g)=\mathfrak{L}_{\lambda^\prime}(g)$ for each $g\in C_c^\infty(Y)$ so that $\mathfrak{L}_{\lambda}=\mathfrak{L}_{\lambda^\prime}$ and hence that $\lambda=\lambda^\prime$.

For b) we see for any $g\in C_c^\infty(Y)$ that
$$\mathcal{P}^h_g(x)=\int_{\widehat{Y}}h(x,\psi)\mathcal{F}_Y(g)(\psi)d\nu_Y$$
which is locally constant since $\psi\mapsto h(x,\psi)\mathcal{F}_Y(g)(\psi)$ is locally constant and compactly supported on $\widehat{Y}$.

For c) we note that $\widehat{Y}$ is discrete and that $\psi\in C_c^\infty(Y)$ with $\mathcal{F}_Y(\psi)=\mathbf{1}_{\left\{\psi\right\}}$ for all $\psi\in\widehat{Y}$. It follows that for all $f\in C_c^\infty(X)$ and all $\psi\in\widehat{Y}$ that
\begin{align*}
\left<\mathfrak{L}_\lambda(\psi),f\right>=&\int_X\mathcal{P}_\psi^h(x)f(a_G)d\mu_X\\
=&\int_X\int_{\widehat{Y}}h(x,\psi^\prime)\mathcal{F}_Y(\psi)(\psi)d\nu_Yf(a_G)d\mu_X\\
=&\int_Xh(x,\psi)f(a_G)d\mu_X\\
=&\int_X\chi_\psi(a_G)f(a_G)d\mu_X
\end{align*}
as required.
\end{proof}

In light of the above, for any $h\in C^\infty(X\times\widehat{Y})$ we henceforth denote by $\mathcal{P}^h$ the smooth family associated to $h$ by Proposition \ref{inversionprop} b) and we denote by $\lambda^h$ the element of $\mathcal{D}(X\times Y)$ associated to $\mathcal{P}^h$ by Proposition \ref{inversionprop} a).


\subsection{Orbital Integtals}\label{OIsubsec}


Following the conventions of \cite{FLN} for a reductive group $J$ over $F$ and $h\in C_c^\infty(J(F))$ we define $\text{Orb}_J(h):J(F)\to\mathbb{C}$ via
$$\text{Orb}_J(h)(\gamma)=\int_{C_J(\gamma)\backslash J(F)}h(\gamma^j)\frac{d\mu_J}{d\mu_{C_J(\gamma)}}$$
for the quotient measure $\frac{d\mu_J}{d\mu_{C_J(\gamma)}}$ as normalized in \cite{FLN} and which is known to converge for all $\gamma\in J(F)$ by \cite{Rao}. 

For $h\in C_c^\infty(G(F))$ we define $\text{St.Orb}_G(h):\mathfrak{A}_G\to\mathbb{C}$ via
$$\text{St.Orb}_G(h)(a_G)=\sum_{\gamma:\pi_G(\gamma)=a_G}\text{Orb}_G(h)(\gamma)$$
so that under our normalization of measures we have $(\pi_G)_\ast(hd\mu_G)=\text{St.Orb}(h)d\eta_G$. We remark that we will not employ the often-used notation $h_G=\text{Orb}_G(h)$ and $h^G=\text{St.Orb}_G(h)$ since we will be working primarily with functions on $\mathfrak{A}_G(F)$ itself and, in our main proofs, will have no need to reference functions on the group $G(F)$ itself. Indeed, we have the following fact which allows us to characterize the functions on $\mathfrak{A}_G(F)^{\text{rss}}$ which are orbital integrals sufficiently for our purposes. We define $C_b^\infty(\mathfrak{A}_G)\subset C^\infty(\mathfrak{A}_G)$ as follows; letting $\iota:C^\infty(\mathfrak{A}_G(F)^{\text{rss}})\to C^\infty(\mathfrak{A}_G(F))$ be the extension by zero map
$$\iota(f)(a_G)=\left\{\begin{array}{ll}
f(a_G)&a_G\in\mathfrak{A}_G(F)^{\text{rss}}\\\\
0&\text{else}
\end{array}\right.$$
we define
$$C_b^\infty(\mathfrak{A}_G)=\{f\in C^\infty(\mathfrak{A}_G^{\text{rss}}(F)):\iota(f)\text{ has compact support}\}.$$
For $f\in C^\infty(\mathfrak{A}_G(F))$ we write abuse notation and write $f\in C_b^\infty(\mathfrak{A}_G)$ if $f|_{\mathfrak{A}_G^{\text{rss}}}\in C_b^\infty(\mathfrak{A}_G)$ so that under this identification we have $C_c^\infty(\mathfrak{A}_G)\subset C_b^\infty(\mathfrak{A}_G)\subset C^\infty(\mathfrak{A}_G)$. We now observe that stable orbital integrals lie in the space $C_b^\infty(\mathfrak{A}_G)$.

\begin{prop}\label{OISL}
For $h\in C_c^\infty(H)$ we have that $\theta_h=|D_G|^\frac{1}{2}\text{\emph{St.Orb}}(h)\in C_b^\infty(\mathfrak{A}_G)$.
\end{prop}
\begin{proof}
To begin, since there are finitely many stable conjugacy classes of maximal tori $S\subset G$ it suffices to show that for all $h\in C_c^\infty(G(F))$ that the function $h^{S,\prime}:S(F)\to \pi_S(S(F))$ defined via $h^{S,\prime}(s)=\text{St.Orb}_G(h)(\pi(s))$ lies in $C_b^\infty(\mathfrak{A}_G)$. In turn, and identifying $\pi_S(S(F))\simeq S(F)/W_S$, it suffices to show that the function $h^S:S(F)^{\text{rss}}\to S(F)^{\text{rss}}/W_S$ defined via $h^S(s)=\text{Orb}_G(h)(s)$ is such that $|D_G|^\frac{1}{2}h^S$ is locally constant with support whose closure in $S(F)/W_S$ is compact.

Denote by $\pi_{S/W_S}:{^{G(F)}S(F)}\to S(F)/W_S$ be the quotient map $\pi_{S/W_S}({^\gamma s})=sW_S$.  From \cite{FLN} we have that there are surjective, open maps $\alpha:S(F)\times S(F)\backslash G(F)\to S(F)^{G(F)}$ given by $(s,S(F)\gamma)\mapsto s^\gamma$ and $\beta:S(F)\times G(F)/S(F)\to S(F)/W_S$ given by $\beta((s,\gamma S(F)))\mapsto sW_s$ which satisfy the property that, for some non-zero constant $C\in\mathbb{C}^\times$ depending on the normalization of measures, we have
$$\beta_\ast\left(fd\mu_S\times\frac{d\mu_G}{d\mu_S}\right)=C\cdot(\pi_{S/W_S})_\ast\left(\alpha_\ast\left(fd\mu_S\times\frac{d\mu_G}{d\mu_S}\right)\right)$$
for all $f\in C_c^\infty(S(F)\times G(F)/S(F))$ as well as that both of the maps $\alpha$ and $\beta$ restricted to $S(F)^{\text{rss}}\times G(F)/S(F)$ are smooth. It follows that $h^S\in C_c^\infty(S(F)^{\text{rss}})$ for $h\in C_c^\infty(G(F)^{\text{rss}})$ and hence that $|D_G|^\frac{1}{2}h^S\in C_c^\infty(S(F)^{\text{rss}})$. 

It remains to show that for all $s\in S(F)\setminus S(F)^{\text{rss}}$ there is a neighbourhood $U_s\subset S(F)$ of $s$ such that $|D_G|^\frac{1}{2}h^S|_{U_s}$ is bounded and locally constant which follows from well-known facts regarding Shalika germ expansions. By \cite{Shalika} and \cite{Vigneras} we have that for $M=C_G(s)^\circ$ there exists a function $h_M\in C_c^\infty(M(F))$ such that $\text{Orb}_M(h_M)|_{S(F)}$ agrees with $\text{Orb}_G(h)|_{S(F)}$ in a neighbourhood of $s$ and such that $\text{Orb}_M(f_M)$ may be written as a linear combination of finitely many Shalika germs $\Gamma_i^M$, $1\le i\le N$ for some positive integer $N$, of $M$. These Shalika germs are locally constant functions and it follows from \cite{KNotes} Theorem 17.9 that the normalized Shalika germs $|D_M|^\frac{1}{2}\Gamma_i^M$ are locally bounded on $S(F)$; the result follows.
\end{proof}

Similarly, for $\widetilde{h}\in C_c^\infty(\widetilde{G}(F))$ we define $\text{Orb}_{\widetilde{G}}(h):G(F)\to\mathbb{C}$ via
$$\text{Orb}_{\widetilde{G}}(\widetilde{h})(\gamma)=\int_{\widetilde{G}(F)/C_{\widetilde{G}(\gamma)}}\widetilde{h}({^g\gamma})\frac{d\mu_{\widetilde{G}}}{d\mu_{C_{\widetilde{G}(\gamma)}}}$$
again normalized as in \cite{FLN} and 
$$\text{St.Orb}_{\widetilde{G}}(\widetilde{h})(a_{\widetilde{G}})=\sum_{\gamma:\pi_{\widetilde{G}}(\gamma)=a_{\widetilde{G}}}\text{Orb}_G(\widetilde{h})(\gamma)$$
so that we similarly have $(\pi_{\widetilde{G}})_\ast(\widetilde{h}d\mu_{\widetilde{G}})=\text{St.Orb}_{\widetilde{G}}(\widetilde{h})d\eta_{\widetilde{G}}$. We define the subspace $C_b^\infty(\mathfrak{A}_{\widetilde{G}})\subset C^\infty(\mathfrak{A}_{\widetilde{G}}(F)^{\text{rss}})$ analogously to the definition of $C_b^\infty(\mathfrak{A}_G)$ above. We now have the following similar result pertaining to $\widetilde{G}(F)$.

\begin{prop}\label{GLproductOI}
For $\theta_{\widetilde{f}}=|D_{\widetilde{G}}|^\frac{1}{2}\text{\emph{st.orb}}(\widetilde{f})|_{\mathfrak{A}_{\widetilde{G}}(F)^{\text{\emph{rss}}}}$ we have that $\theta_{\widetilde{f}}\in C_b^\infty(\mathfrak{A}_{\widetilde{G}})$ for all $\widetilde{f}\in C_c^\infty(\widetilde{G}(F))$. Moreover, for $\widetilde{f}$ with support in $Z_{\widetilde{G}(F)}G(F)_{0^+}$ we have that $\theta_{\widetilde{f}}\in C_c^\infty(Z_{\widetilde{G}}(F))\otimes C_b^\infty(\mathfrak{A}_G)$.
\end{prop}
\begin{proof}
The proof of the first statement is similar to that of Proposition \ref{OISL}. For the second, similarly to the proof of Proposition \ref{OISL}, for each maximal torus $\widetilde{S}\subset\widetilde{G}$ there are maps $\widetilde{\alpha}:\widetilde{S}(F)\times \widetilde{S}(F)\backslash \widetilde{G}(F)\to \widetilde{S}(F)^{\widetilde{G}(F)}$ given by $(s,\widetilde{S}(F)\gamma)\mapsto s^\gamma$ and $\widetilde{\beta}:\widetilde{S}(F)\times \widetilde{S}(F)\backslash \widetilde{G}(F)\to \pi_{\widetilde{S}}(\widetilde{S}(F))$ given by $\widetilde{\beta}((s,\gamma S(F)))\mapsto \pi_{\widetilde{S}}(s)$ satisfying 
$$\widetilde{\beta}_\ast\left(fd\mu_{\widetilde{S}}\times\frac{d\mu_{\widetilde{G}}}{d\mu_{\widetilde{S}}}\right)=\widetilde{C}\cdot(\pi_{\widetilde{S}})_\ast\left(\widetilde{\alpha}_\ast\left(fd\mu_{\widetilde{S}}\times\frac{d\mu_{\widetilde{G}}}{d\mu_{\widetilde{S}}}\right)\right)$$
for some non-zero constant $\widetilde{C}\in\mathbb{C}^\times$. If $f\in C_c^\infty(\widetilde{G}(F))$ has support in $Z_{\widetilde{G}(F)}G(F)_{0^+}$ we have that $f|_{\widetilde{S}(F)^{\widetilde{G}(F)}}$ has support in $Z_{\widetilde{G}(F)}S(F)_{0^+}$ where the map $Z_{\widetilde{G}(F)}\otimes S(F)_{0^+}\to Z_{\widetilde{G}(F)}S(F)_{0^+}$ is an isomorphism. The result follows since we have $\widetilde{\beta}|_{Z_{\widetilde{G}(F)}S(F)_{0^+}\times \widetilde{S}(F)\backslash \widetilde{G}(F)}$ similarly factorizes via
$$\widetilde{\beta}((zs,\widetilde{S}(F)g))=z\pi_S(s)$$
\end{proof}

\section{The $\text{SL}_\ell(F)$ Case}\label{secSL}

In this section we prove Theorem \ref{MAINSL} which will follow immediately from Propositions \ref{SLroundup} and \ref{SLtransferA} and Theorem \ref{SLMainThm} below. We proceed by first computing the stable transfer factor $\Theta_\phi$ associated to the $L$-embedding $\phi:{^LT}\to{^LG}$ defined in $\S$\ref{param} using Proposition \ref{inversionprop}. We then establish the existence of and compute the transfer operator $\mathfrak{S}_\phi$, first for $f\in C_c^\infty(\mathfrak{A}_G(F)^{\text{rss}})$ and, subsequently, for all $f\in C_b^\infty(\mathfrak{A}_G(F))$. Notably, this answers in the affirmative Langlands' Questions A and B of \cite{SetT} for the $L$-embedding $\phi:{^LT}\to{^LG}$. Our argument to extend the definition of the transfer operator from $C_c^\infty(\mathfrak{A}_G(F)^{\text{rss}})$ to $C_b^\infty(\mathfrak{A}_G(F))$ is similar in spirit to that which appears in \cite{SetT} in the case of $\text{SL}_2$. We remark that a novel feature of the stable transfer factor $\Theta_\phi$ in our case is that, unlike in the case of $\text{SL}_2$, it is not represented by a smooth function.

\subsection{Computation of $\Theta_\phi$}\label{SLfamily}

We compute the smooth family $\mathcal{P}_\phi$ associated to the function $_\phi(a_G,\psi)=\chi_\psi(a_G)$ as in $\S$\ref{distsubsec} which, by Proposition \ref{inversionprop}, will yield our desired distribution $\Theta_\phi$, initially considered as an element of $\text{Hom}(C_c^\infty(\mathfrak{A}_G(F)^{\text{rss}}\otimes C_c^\infty(T(F)),\mathbb{C})$. As seen above in $\S$\ref{charformsubsec}, the formulas for the characters $\chi_\psi$ behave very differently on different stable conjugacy classes of $G(F)$ and, moreover, depend crucially on the depth of elements. To clarify the many computations which must be performed, we break up the space $\mathfrak{A}_G(F)$ into various pieces corresponding to stable conjugacy classes and depth. Let $\mathfrak{A}^+=\pi_G(Z_G(F)G(F)_{0^+}\cap G(F)^{\text{rss}})$, $\mathfrak{A}^T=\pi_G(T(F)^{\widetilde{G}(F)})$ and $\mathfrak{A}^A=\pi_G(A(F)^{\widetilde{G}(F)})$, as well as $\mathfrak{A}^{T,-}=\mathfrak{A}^T\setminus \mathfrak{A}^+$ and $\mathfrak{A}^{A,-}=\mathfrak{A}^A\setminus \mathfrak{A}^+$.

We define a function $\Theta^+(a_G,t)$ on $\mathfrak{A}_G\times T(F)$. Despite its seemingly complicated definition we remark that it depends only on the stable conjugacy class of (any lift of) $a_G$, various discriminant functions, as well as the maximal depth of the elements $a_G$ and $t$; we remark that it has close relationship with the Local Character Expansion. Moreover, the function behaves slightly differently and has additional symmetry when $a_G\in\mathfrak{A}^T$. 

We define
\begin{align*}
\Theta^+(a_G,t)=&\left\{\begin{array}{ll}
{\displaystyle\sum_{\mathcal{O}\le\mathcal{O}_\gamma}(A_\mathcal{O}+E_\mathcal{O}q^{\frac{|\Phi_\mathcal{O}|+2\ell-2}{2}(\lceil d^+(a_G)\rceil-1)})B^\mathcal{O}(a_G)}&d^+(a_G)\le d^+(t)<0, \\
&e(a_G,t)=1, a_G\in\mathfrak{A}^+\setminus \mathfrak{A}^T\\\\
{\displaystyle\sum_{\mathcal{O}\le\mathcal{O}_\gamma}(A_\mathcal{O}+D_\mathcal{O}q^{\frac{|\Phi_\mathcal{O}|+2\ell-2}{2}d^+(t)})B^\mathcal{O}(a_G)}&d^+(t)<d^+(a_G)<0, \\
&e(a_G,t)=1, a_G\in\mathfrak{A}^+\setminus \mathfrak{A}^T\\\\
A_{\{0\}}+D_{\{0\}}q^{\frac{\ell^2+\ell-2}{2}\text{min}\{d^+(a_G),d^+(t)\}}&d^+(a_G), d^+(t)<0, \\
&e(a_G,t)=1, a_G\in\mathfrak{A}^T\\\\
{\displaystyle \sum_{\mathcal{O}\le\mathcal{O}_\gamma}(-C_\mathcal{O})B^\mathcal{O}(a_G)}&d^+(a_G)>0, \\
&e(a_G,t)=0, a_G\in\mathfrak{A}^+\setminus \mathfrak{A}^T\\\\
{\displaystyle-\ell+\ell q^\frac{\ell^2-\ell}{2}}&d^+(a_G)>0, \\
&e(a_G,t)=0, a_G\in\mathfrak{A}^T\\\\
0&\text{else}
\end{array}\right.
\end{align*}
for $B^\mathcal{O}$ as defined in $\S$\ref{charformsubsec} and
$$A_\mathcal{O}=\left(|T(\mathfrak{f})|\frac{q^{\frac{|\Phi_\mathcal{O}|}{2}}-1}{q^{\frac{|\Phi_\mathcal{O}|+2\ell-2}{2}}-1}-1\right)C_\mathcal{O}$$
as well as
$$E_\mathcal{O}=|T(\mathfrak{f})|\frac{q^{\frac{|\Phi_\mathcal{O}|+2\ell-2}{2}}-q^{\frac{|\Phi_\mathcal{O}|}{2}}}{q^{\frac{|\Phi_\mathcal{O}|+2\ell-2}{2}}-1}C_\mathcal{O}$$
and
$$D_\mathcal{O}=-|T(\mathfrak{f})|\frac{q^{\frac{|\Phi_\mathcal{O}|}{2}}-1}{q^{\frac{|\Phi_\mathcal{O}|+2\ell-2}{2}}-1}C_\mathcal{O}.$$

We now break our computation into the cases $a_G\notin\mathfrak{A}^T$ and $a_G\in\mathfrak{A}^T$.
\begin{prop}\label{SLnotT}
For $a_G\notin\mathfrak{A}^T$ we have
$$\mathcal{P}(a_G)=\left(\chi_1(a_G)+\Theta^+(a_G,\bullet)\right)d\mu_T.$$
\end{prop}

\begin{proof}
First we suppose $a_G\in\mathfrak{A}^+\setminus \mathfrak{A}^T$. We will show that $\mathcal{F}_{T(F)}(Y(a_G,\bullet))=\chi_1(a_G)+\Theta^+(a_G,\bullet)$ from which it follows that
$$\mathcal{P}(a_G)=\mathcal{F}_{T(F)}(Y(a_G,\bullet)d\nu_T)=\mathcal{F}_{T(F)}(Y(a_G,\bullet))d\mu_T=\left(\chi_1(a_G)+\Theta^+(a_G,\bullet)\right)d\mu_T$$ 
since $\psi\mapsto Y(a_G,\psi)$ is compactly supported by Proposition \ref{CharacterTable}.

Letting $\mathcal{C}_0=\left\{\psi\in\widehat{T(F)}:d(\psi)\le 0\right\}$ and $\mathcal{C}_+=\left\{\psi\in\widehat{T(F)}:d(\psi)> 0\right\}$ we may write
$$\mathcal{F}_{T(F)}(Y(a_G,\bullet))(t)=\int_{\mathcal{C}_0}\chi_\psi(a_G)\psi(t^{-1})d\eta_{\widehat{T}}+\int_{\mathcal{C}_+}\chi_\psi(a_G)\psi(t^{-1})d\eta_{\widehat{T}}$$
and compute these two integrals directly. By Corollary \ref{slstuff} we may choose $\gamma\in G(F)_{0^+}$ and $\zeta\in Z_G(F)$ such that $a_G=\pi_G(\zeta\gamma)$. By Proposition \ref{CharacterTable} we have that $\chi_\psi(\gamma)$ depends only on $d(\psi)$ and throughout this proof we will write $\chi_d(\gamma)=\chi_\psi(\gamma)$ for any $\psi\in\widehat{T(F)}$ with $d(\psi)=d$.

For the integral over $\mathcal{C}_0$, if $e(a_G,t)=0$ we have by Proposition \ref{CharacterTable} that $\chi_{\psi_1}(\gamma)=\chi_{\psi_2}(\gamma)$ for all $\psi_1,\psi_2\in\widehat{T(F)}$ with $d(\psi_1)=d(\psi_2)=0$ so that we have
\begin{align*}
&\int_{\mathcal{C}_0}\chi_\psi(a_G)\psi(t^{-1})d\psi\\
=&\chi_1(\zeta\gamma)+\sum_{\psi\in\widehat{T(F)}:d(\psi)=0}\chi_\psi(\zeta\gamma)\psi(t^{-1})\\
=&\chi_1(\zeta\gamma)+\chi_0(\gamma)\sum_{\psi\in\widehat{T(F)}:d(\psi)=0}\psi(\zeta t^{-1})\\
=&\chi_1(\zeta\gamma)-\chi_0(\gamma)\\
=&\chi_1(a_G)-\sum_{\mathcal{O}\le\mathcal{O}_{a_G}}C_\mathcal{O}\sum_{w\in W_\mathcal{O}(\gamma)}\frac{|D_{M_\mathcal{O}}(\gamma^w)|^\frac{1}{2}}{|D_{\widetilde{G}}(\gamma)|^\frac{1}{2}}.
\end{align*}
On the other hand, if $e(a_G,t)=1$ we use Corollary \ref{slstuff} to write $t=\zeta t^\prime$ and by Proposition \ref{CharacterTable} and Lemma \ref{CharCount} we have
\begin{align*}
&\int_{\mathcal{C}_0}\chi_\psi(a_G)\psi(t^{-1})d\psi\\
=&\chi_1(\zeta\gamma)+\sum_{\psi\in\widehat{T(F)}:d(\psi)=0}\chi_\psi(\zeta\gamma)\psi(t^{-1})\\
=&\chi_1(\zeta\gamma)+\chi^0(\gamma)\sum_{\psi\in\widehat{T(F)}:d(\psi)=0}\psi((t^\prime)^{-1})\\
=&\chi_1(\zeta\gamma)+\chi^0(\gamma)(|T(\mathfrak{f})|-1)\\
=&\chi_1(a_G)+\sum_{\mathcal{O}\le\mathcal{O}_{a_G}}(|T(\mathfrak{f})|-1)C_\mathcal{O}\sum_{w\in W_\mathcal{O}(\gamma)}\frac{|D_{M_\mathcal{O}}(\gamma^w)|^\frac{1}{2}}{|D_{\widetilde{G}}(\gamma)|^\frac{1}{2}}.
\end{align*}

We now consider the integral over $\mathcal{C}_+$. We note in all cases that by Proposition \ref{CharacterTable} and the fact that for any $\gamma^\prime\in T(F)$ with $d(\gamma^\prime)<d$ we have
\begin{equation}\label{VanishFact}\sum_{\psi\in\widehat{T(F)}:d(\psi)=d}\psi(\gamma^\prime)=-\sum_{\psi\in\widehat{T(F)}:d(\psi)<d}\psi(\gamma^\prime)=0\end{equation}
from which it follows that
\begin{equation}\label{BigVanishFact}\int_{\mathcal{C}_+}\chi_\psi(a_G)\psi(t^{-1})d\eta_{\widehat{T}}=\sum_{d=1}^{\lceil d^+(\gamma)\rceil}\sum_{\psi\in\widehat{T(F)}:d(\psi)=d}\chi_\psi(a_G)\psi(t^{-1})\end{equation}
If $e(a_G,t)=0$ we have by (\ref{VanishFact}), (\ref{BigVanishFact}) and the fact that $\chi_\psi(\gamma)$ depends only on $d(\psi)$ by Proposition \ref{CharacterTable} that
\begin{align*}
&\int_{\mathcal{C}_+}\chi_\psi(a_G)\psi(t^{-1})d\psi\\
=&\sum_{d=1}^{\lceil d^+(\gamma)\rceil}\sum_{\psi\in\widehat{T(F)}:d(\psi)=d}\chi_\psi(\zeta\gamma)\psi(t^{-1})\\
=&\sum_{d=1}^{\lceil d^+(\gamma)\rceil}\chi_d(\gamma)\sum_{\psi\in\widehat{T(F)}:d(\psi)=d}\psi(\zeta t^{-1})\\
=&0
\end{align*}
which, combined with our computation of the integral over $\mathcal{C}_0$ above, gives the desired result in this case. 

Now we assume $e(a_G,t)=1$ and write $t=\zeta t^\prime$. If $d^+(a_G)\le d^+(t)$ we have by Proposition \ref{CharacterTable}, (\ref{BigVanishFact}) and Lemma \ref{CharCount} that
\begin{align*}
&\int_{\mathcal{C}_+}\chi_\psi(a_G)\psi(t^{-1})d\psi\\
=&\sum_{\mathcal{O}\le\mathcal{O}_{a_G}}C_\mathcal{O}\sum_{w\in W_\mathcal{O}(\gamma)}\frac{|D_{M_\mathcal{O}}(\gamma^w)|^\frac{1}{2}}{|D_{\widetilde{G}}(\gamma)|^\frac{1}{2}}\sum_{d=1}^{\lceil d^+(a_G)\rceil-1} q^{\frac{|\Phi_\mathcal{O}|}{2}d}\sum_{\psi\in\widehat{T(F)}:d(\psi)=d}\psi((t^\prime)^{-1})\\
=&\sum_{\mathcal{O}\le\mathcal{O}_{a_G}}|T(\mathfrak{f})|(q^{\ell-1}-1)q^{\frac{|\Phi_\mathcal{O}|}{2}}C_\mathcal{O}\sum_{w\in W_\mathcal{O}(\gamma)}\frac{|D_{M_\mathcal{O}}(\gamma^w)|^\frac{1}{2}}{|D_{\widetilde{G}}(\gamma)|^\frac{1}{2}}\sum_{d=0}^{\lceil d^+(a_G)\rceil-2} q^{\frac{|\Phi_\mathcal{O}|+2\ell-2}{2}d}\\
=&\sum_{\mathcal{O}\le\mathcal{O}_{a_G}}|T(\mathfrak{f})|(q^{\ell-1}-1)q^{\frac{|\Phi_\mathcal{O}|}{2}}C_\mathcal{O}\sum_{w\in W_\mathcal{O}(\gamma)}\frac{|D_{M_\mathcal{O}}(\gamma^w)|^\frac{1}{2}}{|D_{\widetilde{G}}(\gamma)|^\frac{1}{2}}\frac{q^{\frac{|\Phi_\mathcal{O}|+2\ell-2}{2}(\lceil d^+(a_G)\rceil-1)}-1}{q^{\frac{|\Phi_\mathcal{O}|+2\ell-2}{2}}-1}\\
=&\sum_{\mathcal{O}\le\mathcal{O}_{a_G}}(-E_\mathcal{O}+E_\mathcal{O}q^{\frac{|\Phi_\mathcal{O}|+2\ell-2}{2}(\lceil d^+(a_G)\rceil-1)})\sum_{w\in W_\mathcal{O}(\gamma)}\frac{|D_{M_\mathcal{O}}(\gamma^w)|^\frac{1}{2}}{|D_{\widetilde{G}}(\gamma)|^\frac{1}{2}}
\end{align*}
which is what we desire given our computation above of the integral over $\mathcal{C}_0$ and since $A_\mathcal{O}=-E_\mathcal{O}+(|T(\mathfrak{f})|-1)C_\mathcal{O}$. Similarly, If $d^+(a_G)> d^+(t)$ we again have by Proposition \ref{CharacterTable}, (\ref{BigVanishFact}) and Lemma \ref{CharCount} and proceeding similarly to the previous case we compute
\begin{align*}
&\int_{\mathcal{C}_+}\chi_\psi(a_G)\psi(t^{-1})d\psi\\
=&\sum_{\mathcal{O}\le\mathcal{O}_{a_G}}C_\mathcal{O}\sum_{w\in W_\mathcal{O}(\gamma)}\frac{|D_{M_\mathcal{O}}(\gamma^w)|^\frac{1}{2}}{|D_{\widetilde{G}}(\gamma)|^\frac{1}{2}}\sum_{d=1}^{d^+(t)-1} q^{\frac{|\Phi_\mathcal{O}|}{2}d}\sum_{\psi\in\widehat{T(F)}:d(\psi)=d}\psi((t^\prime)^{-1})\\
&+\sum_{\mathcal{O}\le\mathcal{O}_{a_G}}C_\mathcal{O}\sum_{w\in W_\mathcal{O}(\gamma)}\frac{|D_{M_\mathcal{O}}(\gamma^w)|^\frac{1}{2}}{|D_{\widetilde{G}}(\gamma)|^\frac{1}{2}}\sum_{d=d^+(t)}^{\lceil d^+(a_G)\rceil-1} q^{\frac{|\Phi_\mathcal{O}|}{2}d}\sum_{\psi\in\widehat{T(F)}:d(\psi)=d}\psi((t^\prime)^{-1})\\
=&\sum_{\mathcal{O}\le\mathcal{O}_{a_G}}(-E_\mathcal{O}+E_\mathcal{O}q^{\frac{|\Phi_\mathcal{O}|+2\ell-2}{2}(d^+(t)-1)})\sum_{w\in W_\mathcal{O}(\gamma)}\frac{|D_{M_\mathcal{O}}(\gamma^w)|^\frac{1}{2}}{|D_{\widetilde{G}}(\gamma)|^\frac{1}{2}}\\
&-\sum_{\mathcal{O}\le\mathcal{O}_{a_G}}|T(\mathfrak{f})|q^{\frac{|\Phi_\mathcal{O}|}{2}}C_\mathcal{O}\sum_{w\in W_\mathcal{O}(\gamma)}\frac{|D_{M_\mathcal{O}}(\gamma^w)|^\frac{1}{2}}{|D_{\widetilde{G}}(\gamma)|^\frac{1}{2}}q^{\frac{|\Phi_\mathcal{O}|+2\ell-2}{2}(d^+(t)-1)}\\
=&\sum_{\mathcal{O}\le\mathcal{O}_{a_G}}(-E_\mathcal{O}+D_\mathcal{O}q^{\frac{|\Phi_\mathcal{O}|+2\ell-2}{2}(d^+(t))})\sum_{w\in W_\mathcal{O}(\gamma)}\frac{|D_{M_\mathcal{O}}(\gamma^w)|^\frac{1}{2}}{|D_{\widetilde{G}}(\gamma)|^\frac{1}{2}}
\end{align*}
as required.

Finally, if $a_G\notin\mathfrak{A}^+\cup\mathfrak{A}^T$ we have by Proposition \ref{CharacterTable} and the definition of $\Theta^+$ that
$$\mathcal{P}_\phi(a_G)=\chi_1(a_G)d\mu_T=\left(\chi_1(a_G)+\Theta^+(a_G,\bullet)\right)d\mu_T$$
since $\chi_\psi(a_G)=0$ for all $\psi\neq 1$ and $\Theta^+(a_G,t)=0$ for all $t\in T(F)$.
\end{proof}
We remark that the map $(a_G,t)\mapsto\chi_1(a_G)$ is smooth and compactly supported and could, in principal, have been incorporated into the definition of the function $\Theta^+$; in addition to being not entirely supported within $\mathfrak{A}^+$, the analogous term in the case of $\widetilde{G}(F)$ must be considered separately; in part this is due to the fact that the determinant of an element of the split torus of $\widetilde{G}(F)$ need not be an $\ell^{\text{th}}$ power. It is for this reason of compatibility between the $G(F)$ and $\widetilde{G}(F)$ cases that we notationally separate this term.

We now perform a similar computation in the case $a_G\in\mathfrak{A}^T$. We define the distribution
$$\delta^T(a_G)=\frac{1}{|D_G(a_G)|^\frac{1}{2}}\sum_{t^\prime\in T(F):\pi(t^\prime)=a_G}\delta_{t^\prime}-\frac{1}{|D_G(a_G)|^\frac{1}{2}}1d\mu_T$$
which is to appear in our computation. We have the following. 
\begin{prop}\label{SLinT}
For $a_G\in\mathfrak{A}^T$ we have
$$\mathcal{P}(a_G)=\delta(a_G)+\Theta^+(a_G,\bullet)d\mu_T.$$
\end{prop}
\begin{proof}
Let $a_G\in\mathfrak{A}^T$ and set $\mathcal{C}_{a_G}=\left\{\psi\in\widehat{T(F)}:0\le d(\psi)<d^+(a_G)\right\}$. We have by Proposition \ref{CharacterTable} that
\begin{align*}
\left<\mathcal{P}(a_G),g\right>=&\left<Y(a_G,\bullet),\mathcal{F}_{T(F)}(g)\right>\\
=&\frac{1}{|D_G(a_G)|^\frac{1}{2}}\sum_{t^\prime:\pi(t^\prime)=a_G}\left<t^\prime,\mathcal{F}_{T(F)}(g)\right>-\frac{1}{|D_G(a_G)|^\frac{1}{2}}\sum_{t^\prime:\pi(t^\prime)=a_G}\left<t^\prime\mathbf{1}_{\{1\}},\mathcal{F}_{T(F)}(g)\right>\\
&+\left<\left(Y(a_G,\bullet)-\frac{1}{|D_G(a_G)|^\frac{1}{2}}\sum_{t^\prime:\pi(t^\prime)=a_G}t^\prime\right)\mathbf{1}_{\mathcal{C}_{a_G}},\mathcal{F}_{T(F)}(g)\right>\\
=&<new,g>+\left<\mathcal{F}_{T(F)}^\vee\left(\left(Y(a_G,\bullet)-\frac{1}{|D_G(a_G)|^\frac{1}{2}}\sum_{t^\prime:\pi(t^\prime)=a_G}t^\prime\right)\mathbf{1}_{\mathcal{C}_{a_G}}\right),g\right>
\end{align*}
so all that remains is to compute is the function $\mathcal{F}_{\widehat{T(F)}}\left(\left(Y(a_G,\bullet)-\frac{1}{|D_G(a_G)|^\frac{1}{2}}\sum_{t^\prime:\pi(t^\prime)=a_G}t^\prime\right)\mathbf{1}_{\mathcal{C}_{a_G}}\right)$. 

If $e(a_G,t)=0$ we have that
$$\mathcal{F}_{\widehat{T(F)}}\left(\left(Y(a_G,\bullet)-\frac{1}{|D_G(a_G)|^\frac{1}{2}}\sum_{t^\prime:\pi(t^\prime)=a_G}t^\prime\right)\mathbf{1}_{\mathcal{C}_{a_G}}\right)(t)=\left\{\begin{array}{ll}
0&d^+(a_G)=0\\\\
\ell(q^\frac{\ell^2-\ell}{2}-1)&d^+(a_G)=0
\end{array}\right.$$
very similarly to how the analogous computation proceeds in the proof of Proposition \ref{SLnotT}.

Now suppose $e(a_G,t)=1$ and let $\gamma\in T(F)$ be such that $\pi_G(\gamma)=a_G$ with $\gamma=\zeta\gamma_+$ for $\zeta\in Z_G(F)$ and $\gamma_+\in T(F)_0$ and further write $t=\zeta t_+$; this may be done by Corollary \ref{slstuff}. By Lemmas \ref{slelliptic} and \ref{CharCount} we have
\begin{align}
\nonumber&\mathcal{F}_{\widehat{T(F)}}\left(-\frac{1}{|D_G(a_G)|^\frac{1}{2}}\sum_{t^\prime:\pi_G(t^\prime)=\pi_G(\gamma_+)}\zeta t^\prime\mathbf{1}_{\mathcal{C}_{a_G}}\right)(\zeta t_+)\\
\nonumber=&-\frac{1}{|D_G(a_G)|^\frac{1}{2}}\sum_{t^\prime:\pi_G(t^\prime)=\pi_G(\gamma_+)}\int_{\mathcal{C}_{a_G}}\psi(t^\prime t_+^{-1})\\
=&-|T(\mathfrak{f})|q^{\frac{\ell^2-\ell}{2}d^+(a_G)}q^{(\ell-1)(d^+(a_G)-1)}\sum_{t^\prime:\pi_G(t^\prime)=\pi_G(\gamma_+)}\mathbf{1}_{t^\prime T(F)_{d^+(a_G)}}(t_+)\label{discrcomp}
\end{align}
We now break our computation of $\mathcal{F}_{\widehat{T(F)}}\left(Y(a_G,\bullet)\mathbf{1}_{\mathcal{C}_{a_G}}\right)$ into cases.

If $d^+(t)<d^+(a_G)$ we have that the contribution of (\ref{discrcomp}) is zero and by Proposition \ref{CharacterTable} and Lemma \ref{CharCount} we have
\begin{align*}
&\int_{\mathcal{C}_{a_G}}\chi_\psi(a_G)\psi(t^{-1})d\eta_{\widehat{T}}\\
=&\sum_{\psi:d(\psi)=0}\chi_\psi(\gamma_0)\psi(t_0^{-1})+\sum_{d=1}^{d^+(t)-1}\sum_{\psi:d(\psi)=d}\chi_\psi(\gamma_0)\psi(t_0^{-1})+\sum_{d=d^+(t)}^{d^+(a_G)-1}\sum_{\psi:d(\psi)=d}\chi_\psi(\gamma_0)\psi(t_0^{-1})\\
=&\ell\left(\sum_{\psi:d(\psi)=0}\psi(t_0^{-1})+\sum_{d=1}^{d^+(t)-1}q^{\frac{\ell^2-\ell}{2}d}\sum_{\psi:d(\psi)=d}\psi(t_0^{-1})+\sum_{d=d^+(t)}^{d^+(a_G)-1}q^{\frac{\ell^2-\ell}{2}d}\sum_{\psi:d(\psi)=d}\psi(t_0^{-1})\right)\\
=&\ell\left(|T(\mathfrak{f})|-1+|T(\mathfrak{f})|q^{\frac{\ell^2-\ell}{2}}(q^{\ell-1}-1)\sum_{d=0}^{d^+(t)-2}q^{\frac{\ell^2+\ell-2}{2}d}-|T(\mathfrak{f})|q^{\frac{\ell^2-\ell}{2}}q^{\frac{\ell^2+\ell-2}{2}(d^+(t)-1)}\right)\\
=&\ell\left(|T(\mathfrak{f})|\frac{q^{\frac{\ell^2-\ell}{2}}-1}{q^\frac{\ell^2+\ell-2}{2}-1}-1\right)+\ell|T(\mathfrak{f})|q^{\frac{\ell^2-\ell}{2}}\left(\frac{q^{\ell-1}-1}{q^\frac{\ell^2+\ell-2}{2}-1}-1\right)q^{\frac{\ell^2+\ell-2}{2}(d^+(t)-1)}\\
=&\ell\left(|T(\mathfrak{f})|\frac{q^{\frac{\ell^2-\ell}{2}}-1}{q^\frac{\ell^2+\ell-2}{2}-1}-1\right)-\ell|T(\mathfrak{f})|\frac{q^\frac{\ell^2+\ell-2}{2}-1}{q^\frac{\ell^2+\ell-2}{2}-1}q^{\frac{\ell^2+\ell-2}{2}d^+(t)}
\end{align*}
as required. 

If $d^+(t)\ge d^+(a_G)$ we instead have
\begin{align*}
&\int_{\mathcal{C}_{a_G}}\chi_\psi(a_G)\psi(t^{-1})d\eta_{\widehat{T}}\\
=&\sum_{\psi:d(\psi)=0}\chi_\psi(\gamma_0)\psi(t_0^{-1})+\sum_{d=1}^{d^+(a_G)-1}\sum_{\psi:d(\psi)=d}\chi_\psi(\gamma_0)\psi(t_0^{-1})\\
=&\ell\left(\sum_{\psi:d(\psi)=0}\psi(t_0^{-1})+\sum_{d=1}^{d^+(a_G)-1}q^{\frac{\ell^2-\ell}{2}d}\sum_{\psi:d(\psi)=d}\psi(t_0^{-1})\right)\\
=&\ell\left(|T(\mathfrak{f})|-1+|T(\mathfrak{f})|q^{\frac{\ell^2-\ell}{2}}(q^{\ell-1}-1)\sum_{d=0}^{d^+(t)-2}q^{\frac{\ell^2+\ell-2}{2}d}\right)\\
=&\ell\left(|T(\mathfrak{f})|\frac{q^{\frac{\ell^2-\ell}{2}}-1}{q^\frac{\ell^2+\ell-2}{2}-1}-1\right)+\ell|T(\mathfrak{f})|q^{\frac{\ell^2-\ell}{2}}\frac{q^{\ell-1}-1}{q^\frac{\ell^2+\ell-2}{2}-1}q^{\frac{\ell^2+\ell-2}{2}(d^+(t)-1)}\\
\end{align*}
and we have that (\ref{discrcomp}) contributes the term $-|T(\mathfrak{f})|q^{\frac{\ell^2-\ell}{2}}q^{\frac{\ell^2+\ell-2}{2}(d^+(t)-1)}$ so that we are done.
\end{proof}

We now encapsulate the computations performed in $\S$\ref{SLfamily} as follows.
\begin{prop}\label{SLroundup}
For $\Theta_\phi=\lambda_{\mathcal{P}_\phi}$ we have that $\mathfrak{L}_{\Theta_\phi}(\psi)=\chi_\psi d\eta_G$ for all $\chi\in\widehat{T(F)}$. Moreover we have that $\Theta_\phi=\Theta_\phi^++\Theta_\phi^A+\Theta_\phi^T$ for 
$$\Theta_\phi^+=\Theta^+ d\eta_G\times d\mu_T$$
for $\Theta^+$ defined in $\S$\ref{SLfamily} as well as
$$\Theta_\phi^A=\chi_1d\eta_G\otimes 1d\mu_T$$
and
$$\left<\Theta_\phi^T,f\otimes g\right>=\sum_{\sigma\in\Gamma_{E|F}}\int_{T(F)}\pi_T^\ast(f)(t)g(t^\sigma)dt$$
for all $f\in C_c^\infty(\mathfrak{A}_G(F)^{\text{\emph{rss}}})$ and $g\in C_c^\infty(T(F))$. 
\end{prop}
\begin{proof}
This follows from Propositions \ref{SLnotT} and \ref{SLinT} and \ref{inversionprop} c), as well as the fact for $H\in C^\infty(\mathfrak{A}_G(F)^{\text{rss}}\times T(F))$ we have that
$$\mathfrak{L}_{H d\eta_G\times d\mu_T}(g)=\left(\int_{T(F)}H(\bullet,t)g(t)d\mu_T\right)d\eta_G$$
and that it follows easily by direct computation that the distribution $\Theta^T$ is that which corresponds to the smooth family $a_G\mapsto \delta^T(a_G)$ as defined in $\S$\ref{SLfamily}.
\end{proof}

\subsection{Stable Transfer for $G(F)$}\label{transfersubsecSL}

We now use our computation of the stable transfer factor $\Theta_\phi$ to define a corresponding stable transfer operator. For $f\in C_c^\infty(\mathfrak{A}_G(F)^{\text{rss}})$ the existence of a transfer operator and a computation thereof follows immediately from Proposition \ref{SLroundup}; we remark that the following is the answer to Langlands' Question A of \cite{SetT} for $\phi$ in the affirmative.
\begin{prop}\label{SLtransferA}
There exists a linear map $\mathfrak{S}_\phi^{\text{\emph{rss}}}:C_c^\infty(\mathfrak{A}_G(F)^{\text{\emph{rss}}})\to C_c^\infty(T(F))$ satisfying $\mathfrak{R}_\Theta(f)=\mathfrak{S}_\phi^{\text{\emph{rss}}}(f)d\eta_G$ where we may decompose $\mathfrak{S}_\phi^{\text{\emph{rss}}}=\mathfrak{S}^++\mathfrak{S}^A+\mathfrak{S}^T$ for 
$$\mathfrak{S}^+(f)(t)=\int_{\mathfrak{A}_G(F)}\Theta^+(a_G,t)f(a_G)d\eta_G,$$
$$\mathfrak{S}^A(f)(t)=\mathfrak{S}^A(f)=\int_{\mathfrak{A}_G(F)}\chi_1(a_G)f(a_G)d\eta_G$$
and
$$\mathfrak{S}^T(f)(t)=\pi_T^\ast(f)(t)$$
where $\mathfrak{S}^+$ is supported on $\mathfrak{A}^+$, $\mathfrak{S}^A$ is supported on $\mathfrak{A}^A$, and $\mathfrak{S}^T$ is supported on $\mathfrak{A}^T$.
\end{prop}
\begin{proof}
This follows from Proposition \ref{SLroundup}, the fact that for $H\in C^\infty(\mathfrak{A}_G(F)^{\text{rss}}\times T(F))$ we have that
$$\mathfrak{R}_{H d\eta_G\times d\mu_T}(f)=\left(\int_{\mathfrak{A}_G(F)}H(a_G,\bullet)f(a_H)d\eta_G\right)d\mu_T$$
as well as that direct computation yields $\mathfrak{R}_{\Theta^T}(f)=\pi_T^\ast(f)d\mu_T$.
\end{proof}

It remains to define a transfer operator $\mathfrak{S}_\phi$ which is defined on $C_b^\infty(\mathfrak{A}_G)$ which restricts to $\mathfrak{S}_\phi^{\text{rss}}$ for functions supported on a compact subset of $\mathfrak{A}_G(F)^{\text{rss}}$. With this operator defined, one immediately obtains the extension of $\Theta_\phi$ to $C_b^\infty(\mathfrak{A}_G)\otimes C_c^\infty(T(F))$ by setting
$$\left<\Theta_\phi,f\otimes g\right>=\int_{T(F)}\mathfrak{S}_\phi(f)(t)g(t)d\mu_T$$
so that indeed computing the appropriate transfer operator is all which remains to be done.

For $g\in C^\infty(T(F)^{\text{rss}})$ we define $\iota(g):T(F)\to \mathbb{C}$ to the given by $\iota(g)(t)=g(t)$ for $t\in T(F)^{\text{rss}}$ and $\iota(g)(t)=0$ otherwise. Furthermore, we set
$$C_{\mu_T}^\infty(T(F))=\left\{g\in C^\infty(T(F)^{\text{rss}}):\text{supp}(\iota(g))\text{ is compactly supported, }\int_{T(F)}|\iota(g)(t)|d\mu_T<\infty\right\}$$
where, attaining the first inclusion by restricting functions on $T(F)$ to the open subset $T(F)^{\text{rss}}$, we identify $C_c^\infty(T(F))\subset C_{\mu_T}^\infty(T(F))\subset C^\infty(T(F))$.

We may now show the existence of the operator $\mathfrak{S}_\phi$ and answer in the affirmative Langlands' Question B for $\phi$.
\begin{thm}\label{SLMainThm}
There exists a linear map $\mathfrak{S}_\phi:C_b^\infty(\mathfrak{A}_G(F))\to C_{\mu_T}^\infty(T(F))$ with ${\mathfrak{S}_\phi}|_{C_c^\infty(\mathfrak{A}_G(F)^{\text{\emph{rss}}})}=\mathfrak{S}_\phi^{\text{\emph{rss}}}$ satisfying 
$$\int_{\mathfrak{A}_G(F)}f(a_G)\chi_\psi(a_G)d\eta_G=\int_{T(F)}\mathfrak{S}_\phi(f)(t)\psi(t)dt$$
for each $\psi\in\widehat{T(F)}$.
\end{thm}
\begin{proof}
For any $f\in C_b^\infty(\mathfrak{A}_G(F))$ we have that $f=\sum_{S\in\mathcal{S}}f|_{\pi(S(F)^{\text{rss}})}$ for a set $\mathcal{S}$ of representatives of stable conjugacy classes of maximal tori of $G(F)$; it thus suffices to define $\mathfrak{S}(f_S)$ for $f_S=f|_{\pi(S(F)^{\text{rss}})}$ for each $S\in\mathcal{S}$. Write $f_{S,-1}={f_S}|_{\pi_G(S(F)\setminus Z_G(F)S(F)_{0^+})}$ and for $m\ge 0$ set $f_{S,m}={f_S}|_{\pi(Z_G(F)S(F)_{m^+}\setminus Z_G(F)S(F)_{(m+1)^+}})$ so that $f_S=f_{S,-1}+\sum_{m=0}^\infty f_{S,m}$. Since $f\chi_\psi\in L^1(\mathfrak{A}_G(F),\eta_G)$ we have that
\begin{align*}
\int_{\mathfrak{A}_G(F)}f_S(a_G)\chi_\psi(a_G)d\eta_G=&\int_{\mathfrak{A}_G(F)}{f_{S,-1}}(a_G)\chi_\psi(a_G)d\eta_G+\sum_{m=0}^\infty\int_{\mathfrak{A}_G(F)}f_{S,m}(a_G)\chi_\psi(a_G)d\eta_G\\
=&\int_{T(F)}\mathfrak{S}_\phi^{\text{rss}}(f_{S,-1})(t)\psi(t)d\mu_T+\sum_{m=0}^\infty\int_{T(F)}\mathfrak{S}_\phi^{\text{rss}}(f_{S,m})(t)\psi(t)d\mu_T
\end{align*}
and hence by setting
\begin{equation}\label{Sb}\mathfrak{S}_\phi(f_S)=\mathfrak{S}_\phi^{\text{rss}}(f_{S,-1})+\sum_{m=0}^\infty\mathfrak{S}_\phi^{\text{rss}}(f_{S,m})\end{equation}
it suffices to show that the infinite sum on the right hand side of (\ref{Sb}) converges to an element of $C_{\mu_T}^\infty(T(F))$. We now have various cases depending on the stable class of the torus $S$.

First suppose $S$ is neither stably conjugate to neither $A$ nor $T$ so that $\mathfrak{S}_\phi^{\text{rss}}(f_{S,m})=\mathfrak{S}^+(f_{S,m})$ for all $m\ge -1$ by Proposition \ref{SLroundup} and notably that $\mathfrak{S}_\phi^{\text{rss}}(f_S)(t)=0$ for $t\notin Z_G(F)T(F)_{0^+}$. Denote by $\mathcal{O}_S$ the nilpotent orbit $\mathcal{O}_s$ associated to $s$ for any $s\in S(F)^{\text{rss}}$ and denote by $W_\mathcal{O}(S)$ the set $W_\mathcal{O}(s)$ for any $s\in S(F)^{\text{rss}}$ and $\mathcal{O}\le\mathcal{O}_S$. Moreover, suppose $N_{f,S}$ is an upper bound for $\pi_S^\ast(f_S)$. For $k>0$, $\zeta\in Z_G(F)$ and $t\in \zeta T(F)_k\setminus \zeta T(F)_{k+1}$ we have that
\begin{align*}
\sum_{m=0}^\infty\mathfrak{S}(f_{S,m})(t)=&\sum_{m=0}^\infty\int_{\pi(Z_G(F)S(F)_{m^+}\setminus Z_G(F)S(F)_{(m+1)^+})}f_S(a_G)\Theta^+(a_G,t)d\eta_G\\
=&\sum_{\zeta^\prime\in Z_G(F):\zeta^\prime\neq\zeta}\sum_{m=0}^\infty\int_{\pi(\zeta^\prime S(F)_{m^+}\setminus \zeta^\prime S(F)_{(m+1)^+})}f_S(a_G)\Theta^+(a_G,t)d\eta_G\\
&+\sum_{m=0}^\infty\int_{\pi(\zeta S(F)_{m^+}\setminus \zeta S(F)_{(m+1)^+})}f_S(a_G)\Theta^+(a_G,t)d\eta_G
\end{align*}
where for any $\zeta^\prime\neq\zeta$ we have
\begin{align*}
&\sum_{m=0}^\infty\int_{\pi(\zeta^\prime S(F)_{m^+}\setminus \zeta^\prime S(F)_{(m+1)^+})}f_S(a_G)\Theta^+(a_G,t)d\eta_G\\
=&-\sum_{\mathcal{O}\le\mathcal{O}_S}\sum_{w\in W_\mathcal{O}(S)}C_\mathcal{O} \sum_{m=0}^\infty\int_{\pi(\zeta^\prime S(F)_{m^+}\setminus \zeta^\prime S(F)_{(m+1)^+})}f_S(a_G)\frac{|D_{^wM_\mathcal{O}}(a_G)|^\frac{1}{2}}{|D_G(a_G)|^\frac{1}{2}}d\eta_G\\
=&-\sum_{\mathcal{O}\le\mathcal{O}_S}\sum_{w\in W_\mathcal{O}(S)}C_\mathcal{O} \sum_{m=0}^\infty\int_{\zeta^\prime S(F)_{m^+}\setminus \zeta^\prime S(F)_{(m+1)^+}}\pi_S^\ast(f_S)(s)|D_{^wM_\mathcal{O}}(s)|^\frac{1}{2}d\mu_S
\end{align*}
so that by Lemma \ref{TorusBoundLemma} c) we have
\begin{align*}
&\left|\sum_{m=0}^\infty\int_{\pi(\zeta^\prime S(F)_{m^+}\setminus \zeta^\prime S(F)_{(m+1)^+})}f_S(a_G)\Theta^+(a_G,t)d\eta_G\right|\\
\le&N_{f,S}\sum_{\mathcal{O}\le\mathcal{O}_S}\sum_{w\in W_\mathcal{O}(S)}|C_\mathcal{O}| \sum_{m=0}^\infty\int_{\zeta^\prime S(F)_{m^+}\setminus \zeta^\prime S(F)_{(m+1)^+}}|D_{M_\mathcal{O}}(s^w)|^\frac{1}{2}d\mu_S\\
\le&N_{f,S}Q_S\sum_{\mathcal{O}\le\mathcal{O}_S}|W_\mathcal{O}(S)||C_\mathcal{O}| \sum_{m=0}^\infty q^{-m\left(\frac{|\Phi_\mathcal{O}|}{2}+\ell-1\right)}\\
=&N_{f,S}Q_S\sum_{\mathcal{O}\le\mathcal{O}_S}|W_\mathcal{O}(S)||C_\mathcal{O}| \frac{1}{1-q^{-\left(\frac{|\Phi_\mathcal{O}|}{2}+\ell-1\right)}}.
\end{align*}
We similarly have
\begin{align}
\nonumber&\sum_{m=0}^\infty\int_{\pi(\zeta S(F)_{m^+}\setminus \zeta S(F)_{(m+1)^+})}f_S(a_G)\Theta^+(a_G,t)d\eta_G\\
=&\sum_{\mathcal{O}\le\mathcal{O}_S}\sum_{w\in W_\mathcal{O}(S)}A_\mathcal{O} \sum_{m=0}^\infty\int_{\pi(\zeta^\prime S(F)_{m^+}\setminus \zeta^\prime S(F)_{(m+1)^+})}f_S(a_G)\frac{|D_{^wM_\mathcal{O}}(a_G)|^\frac{1}{2}}{|D_G(a_G)|^\frac{1}{2}}d\eta_G\label{Aline}\\
&+\sum_{\mathcal{O}\le\mathcal{O}_S}\sum_{w\in W_\mathcal{O}(S)}E_\mathcal{O}\sum_{m=0}^{k-1}\int_{\pi(\zeta^\prime S(F)_{m^+}\setminus \zeta^\prime S(F)_{(m+1)^+})}f_S(a_G)q^{\frac{|\Phi_M|+2\ell-2}{2}(\lceil d^+(a_G)\rceil-1)}\frac{|D_{^wM_\mathcal{O}}(a_G)|^\frac{1}{2}}{|D_G(a_G)|^\frac{1}{2}}d\eta_G\label{Eline}\\
&+\sum_{\mathcal{O}\le\mathcal{O}_S}\sum_{w\in W_\mathcal{O}(S)}D_\mathcal{O}q^{\frac{|\Phi_{M_\mathcal{O}}|+2\ell-2}{2}k}\sum_{m=k}^\infty\int_{\pi(\zeta^\prime S(F)_{m^+}\setminus \zeta^\prime S(F)_{(m+1)^+})}f_S(a_G)\frac{|D_{^wM_\mathcal{O}}(a_G)|^\frac{1}{2}}{|D_G(a_G)|^\frac{1}{2}}d\eta_G\label{Dline}
\end{align}
where the convergence of (\ref{Aline}) follows similarly to the argument above, for (\ref{Dline}) we, again by Lemma \ref{TorusBoundLemma} c), have,
\begin{align*}
&\left|\sum_{\mathcal{O}\le\mathcal{O}_S}\sum_{w\in W_\mathcal{O}(S)}D_\mathcal{O}\sum_{m=k}^\infty\int_{\pi(\zeta^\prime S(F)_{m^+}\setminus \zeta^\prime S(F)_{(m+1)^+})}f_S(a_G)\frac{|D_{^wM_\mathcal{O}}(a_G)|^\frac{1}{2}}{|D_G(a_G)|^\frac{1}{2}}d\eta_G\right|\\
\le&N_{f,S}\sum_{\mathcal{O}\le\mathcal{O}_S}\sum_{w\in W_\mathcal{O}(S)}|D_\mathcal{O}|q^{\frac{|\Phi_{M_\mathcal{O}}|+2\ell-2}{2}k}\sum_{m=k}^\infty\int_{\pi(\zeta^\prime S(F)_{m^+}\setminus \zeta^\prime S(F)_{(m+1)^+})}|D_{M_\mathcal{O}}(s^w)|^\frac{1}{2}d\mu_S\\
\le&N_{f,S}Q_S\sum_{\mathcal{O}\le\mathcal{O}_S}\sum_{w\in W_\mathcal{O}(S)}|D_\mathcal{O}|q^{\frac{|\Phi_{M_\mathcal{O}}|+2\ell-2}{2}k}\sum_{m=k}^\infty q^{-m\left(\frac{|\Phi_\mathcal{O}|}{2}+\ell-1\right)}\\
=&N_{f,S}Q_S\sum_{\mathcal{O}\le\mathcal{O}_S}\sum_{w\in W_\mathcal{O}(S)}|D_\mathcal{O}|\frac{1}{1-q^{-\left(\frac{|\Phi_\mathcal{O}|}{2}+\ell-1\right)}}
\end{align*}
and similarly for (\ref{Eline}) we have
\begin{align*}
&\left|\sum_{\mathcal{O}\le\mathcal{O}_S}\sum_{w\in W_\mathcal{O}(S)}E_\mathcal{O}\sum_{m=0}^{k-1}\int_{\pi(\zeta^\prime S(F)_{m^+}\setminus \zeta^\prime S(F)_{(m+1)^+})}f_S(a_G)q^{\frac{|\Phi_M|+2\ell-2}{2}(\lceil d^+(a_G)\rceil-1)}\frac{|D_{^wM_\mathcal{O}}(a_G)|^\frac{1}{2}}{|D_G(a_G)|^\frac{1}{2}}d\eta_G\right|\\
\le&N_{f,S}\sum_{\mathcal{O}\le\mathcal{O}_S}\sum_{w\in W_\mathcal{O}(S)}|E_\mathcal{O}|\sum_{m=0}^{k-1} q^{\frac{|\Phi_{M_\mathcal{O}}|+2\ell-2}{2}d^+(s)}\int_{\pi(\zeta^\prime S(F)_{m^+}\setminus \zeta^\prime S(F)_{(m+1)^+})}|D_{M_\mathcal{O}}(s^w)|^\frac{1}{2}d\mu_S\\
\le&N_{f,S}\sum_{\mathcal{O}\le\mathcal{O}_S}\sum_{w\in W_\mathcal{O}(S)}|E_\mathcal{O}|Q_Sk\\
=&N_{f,S}\sum_{\mathcal{O}\le\mathcal{O}_S}\sum_{w\in W_\mathcal{O}(S)}|E_\mathcal{O}|Q_Sd^+(t)
\end{align*}
It follows that $\mathfrak{S}_b(f_S)$ is integrable since it is bounded above by a function of the form $t\mapsto R_1+R_2d^+(t)$ for constants $R_1,R_2\in\mathbb{C}$. Moreover, it is supported on $Z_G(F)T(F)_{0^+}$ and is locally constant since the above computations show that its values depend only on $d^+(t)$ and the element $\zeta\in Z_G(F)$ with $t\in \zeta T(F)_{0^+}$.

Now suppose $S$ is stably conjugate to $A$; it suffices to suppose $S=A$. The above handles the convergence of $\mathfrak{S}^+(f_A)$ and $\mathfrak{S}^T(f_A)=0$ so it remains to consider the convergence of $\mathfrak{S}^A(f_A)$. We indeed have convergence since
$$\mathfrak{S}^A(f_A)=\int_{\mathfrak{A}_G(F)}\frac{u_1(a_G)}{|D_G(a_G)|^\frac{1}{2}}f_A(a_G)d\eta_G=\int_{A(F)}u_1(a)f_A(a)d\mu_A$$
and $u_1f_A$ is bounded on $A(F)$ with support contained in a compact set.

Finally, for $S$ is stably conjugate to $T$ so that we suppose $S=T$, we similarly need only consider the convergence of $\mathfrak{S}^T(f_T)$ since $\mathfrak{S}^A(f_T)=0$ and the case of $\mathfrak{S}^+(f_T)$ follows as above. We observe that
$$\mathfrak{S}^T(f_T)=\pi_T^\ast(f_T)$$
so that we indeed have $\mathfrak{S}_\phi(f_T)\in C_{\mu_T}^\infty(T(F))$. 
\end{proof}

\section{The $\widetilde{G}=\text{GL}_\ell$ Case}\label{secGL}

In this section we prove Theorem \ref{MAINGL} which follows immediately from Proposition \ref{GLQA} and Theorem \ref{transferforGL}. While our methods are similar to those of $\S$\ref{secSL}, our work is complicated by the fact that $Z_{\widetilde{G}}(F)$ and $\widetilde{T}(F)$ are not compact. Moreover, while we can closely relate the necessary work for $\widetilde{\phi}$ to that done for $\phi$, this process is must be handled carefully given the fact that the multiplication maps $Z_{\widetilde{G}}(F)\times G(F)\to \widetilde{G}(F)$ and $Z_{\widetilde{G}}(F)\times T(F)\to \widetilde{T}(F)$ are neither injective nor surjective. These various difficulties are dealt with in $\S$\ref{helloGL}.

\subsection{Distributions on $\mathfrak{A}_{\widetilde{G}}(F)^{\text{rss}}\times\widetilde{T}(F)$}\label{helloGL}

Similarly to $\S$\ref{secSL} we define the subsets $\mathfrak{A}_{\widetilde{G}}^+=\pi_{\widetilde{G}}(Z_{\widetilde{G}}(F)\widetilde{G}(F)_{0^+})$ and moreover note by Corollary \ref{rootlem} that $Z_{\widetilde{G}}(F)\widetilde{G}(F)_{0^+}=Z_{\widetilde{G}}(F)G(F)_{0^+}$, $\mathfrak{A}_{\widetilde{G}}^T=\pi_{\widetilde{G}}(\widetilde{T}(F))$ and $\mathfrak{A}_{\widetilde{G}}^A=\pi_{\widetilde{G}}(\widetilde{A}(F))$. Also, let $\mathfrak{A}_{\widetilde{G}}(F)^{\det(\widetilde{T})}$ be the subset of elements $a_{\widetilde{G}}\in\mathfrak{A}_{\widetilde{G}}(F)$ with $\det(a_{\widetilde{G}})\in\det(\widetilde{T}(F))$. 

Unlike in the case for $G=\text{SL}_\ell$ where characters of $T(F)$ are compactly supported functions thereupon, for $\widetilde{\psi}\in\widehat{\widetilde{T}(F)}$ we necessarily have $\widetilde{\psi}\notin C_c^\infty(\widetilde{T}(F))$. To deal with this fact, we define a subset $\Delta_{\det}\subset\mathfrak{A}_{\widetilde{G}}(F)\times\widetilde{T}(F)$ via
$$\Delta_{\det}=\left\{(a_{\widetilde{G}},t)\in\mathfrak{A}_{\widetilde{G}}(F)^{\text{rss}}\times\widetilde{T}(F):\det(a_{\widetilde{G}})=\det(t)\right\}$$
and let $\mathcal{D}(\mathfrak{A}_{\widetilde{G}}\times\widetilde{T},\Delta_{\det})\subset\mathcal{D}(\mathfrak{A}_{\widetilde{G}}(F)^{\text{rss}}\times\widetilde{T}(F))$ be the space of distributions supported on $\Delta_{\det}$. We will make fundamental use of the following extension property of $\mathcal{D}(\widetilde{G}\times\widetilde{T},\Delta_{\det})$.

\begin{lem}\label{extension}
Any $\lambda\in\mathcal{D}(\mathfrak{A}_{\widetilde{G}}(F)^{\text{\emph{rss}}}\times\widetilde{T}(F),\Delta_{\det})$ may be extended to a functional on $C_c^\infty(\mathfrak{A}_{\widetilde{G}}(F)^{\text{\emph{rss}}})\otimes C^\infty(\widetilde{T}(F))$ via $\left<\lambda,f\otimes g\right>=\left<\lambda,(f\otimes g)\mathbf{1}_{U_{\Delta_{\det}}}\right>$ for
$$U_{\Delta_{\det}}=\left\{(a_{\widetilde{G}},t)\in\mathfrak{A}_{\widetilde{G}}(F)^{\text{\emph{rss}}}\times\widetilde{T}(F):\det(a_{\widetilde{G}})\in\det(t)\mathcal{O}_F^\times\right\}.$$
Moreover, the corresponding operator $\mathfrak{L}_\lambda:C_c^\infty(\widetilde{T}(F))\to\mathcal{D}(\mathfrak{A}_{\widetilde{G}})$ extends to $C^\infty(\widetilde{T}(F))$ by setting
$$\mathfrak{L}_\lambda(g)=\lim_{m\to\infty}\mathfrak{L}_\lambda\left(g\mathbf{1}_{\left\{t\in\widetilde{T}(F):\text{\emph{ord}}(\det(t))\ge -m\right\}}\right).$$
\end{lem}
\begin{proof}
Let $f\in C_c^\infty(\mathfrak{A}_{\widetilde{G}}(F)^{\text{rss}})$ and $g\in C^\infty(\widetilde{T}(F))$. For $B_f=\left\{t\in \widetilde{T}(F):\det(t)\in\det(\text{supp}(f))\mathcal{O}_F^\times\right\}$ we have that the support of $(f\otimes g)\mathbf{1}_{U_{\Delta_{\det}}}$ is contained in $\text{supp}(f)\times B_f$ and hence is compactly supported. Moreover, there exists an $m_0$ such that $B_f\subset \left\{t\in\widetilde{T}(F):\text{ord}(\det(t))\ge -m\right\}$ and hence that
$$\left<\lambda,f\otimes g\right>=\left<\mathfrak{L}_\lambda\left(g\mathbf{1}_{\left\{t\in\widetilde{T}(F):\text{ord}(\det(t))\ge -m\right\}}\right),f\right>$$
for all $m\ge m_0$.
\end{proof}
We remark that, while the extension defined in Lemma \ref{extension} appears to depend on various choices, it can easily be seen to have the following uniqueness property: if $g\in C^\infty(\widetilde{T}(F))$ is such that $g=\sum_{m=1}^\infty g_m$ for $g_m\in C_c^\infty(\widetilde{T}(F))$ with $\text{supp}(g_{m_1})\cap\text{supp}(g_{m_2})=\emptyset$ for $m_1\neq m_2$ we have that $\left<\lambda,f\otimes g\right>=\sum_{m=1}^\infty\left<\lambda,f\otimes g_m\right>$ for any $f\in C_c^\infty(\mathfrak{A}_{\widetilde{G}}(F)^{\text{rss}})$ and that $\mathfrak{L}_\lambda(g)=\sum_{m=1}^\infty\mathfrak{L}_\lambda(g_m)$.

We now have a number of facts relating to the notions above which will be used in our computations throughout $\S$\ref{GLfamilycomps}. In essence, our goal is to relate various distributions on $\mathfrak{A}_G(F)\times T(F)$ to those on $\mathfrak{A}_{\widetilde{G}}(F)\times\widetilde{T}(F)$ but this must be done carefully given that neither of the multiplication maps $Z_{\widetilde{G}}(F)\times G(F)\to \widetilde{G}(F)$ nor $Z_{\widetilde{G}}(F)\times T(F)\to \widetilde{T}(F)$ are injective or surjective. Even so, we may establish a number of results which will be sufficient for our purposes.

For $z\in Z_{\widetilde{G}}(F)$ we define $\mathcal{M}_z:C_c^\infty(\widetilde{T}(F))\to C_c^\infty(T(F))$ via $\mathcal{M}_z(g)(t)=g(zt)$. Furthermore, for $\lambda\in\mathcal{D}(T(F))$ we define the distribution $\delta_z\otimes\lambda\in\mathcal{D}(\widetilde{T}(F))$ via $\left<\delta_z\otimes\lambda,g\right>=\left<\lambda,\mathcal{M}_z(g)\right>$. We remark that this definition is indeed a minor abuse of notation given that $Z_{\widetilde{G}}(F)\times T(F)\to \widetilde{T}(F)$ is not bijective; even so, this convention shouldn't cause any ambiguity and, moreover, will lessen notational clutter. Also, for $h^Z\in C_c^\infty(Z_{\widetilde{G}}(F))$ and $h^0\in C_c^\infty(T(F))$ define a function $(h^Z\otimes h^0)^{Z_{\widetilde{G}}(F)T(F)}\in C_c^\infty(\widetilde{T}(F))$ supported on $Z_{\widetilde{G}}(F)T(F)$ via
$$(h^Z\otimes h^0)^{Z_{\widetilde{G}}(F)T(F)}(zt)=\sum_{\alpha\in Z_G(F)}h^Z(\alpha z)h^0(\alpha^{-1} t)$$
where we note that $(h^Z\otimes h^0)^{Z_{\widetilde{G}}(F)T(F)}$ has the property that
$$\int_{\widetilde{T}(F)}(h^Z\otimes h^0)^{Z_{\widetilde{G}}(F)T(F)}(y)g(y)d\mu_{\widetilde{T}}=\int_{Z_{\widetilde{G}}(F)}\int_{T(F)}h^Z(z)h^0(t)g(zt)d\mu_Td\mu_{Z_{\widetilde{G}}}$$
for all $g\in C_c^\infty(\widetilde{T}(F))$.



The following Lemma will be used fundamentally in our computation of $\Theta_{\widetilde{\phi}}$. It is through this Lemma that we may view various pieces of $\Theta_{\widetilde{\phi}}$ as natural extensions of analogous pieces of the distribution $\Theta_\phi$.

\begin{lem}\label{GLdistdecomp}
Let $H\in C^\infty(\mathfrak{A}_G(F)\times\widehat{T(F)})$ be supported on $\mathfrak{A}_G^+\times \widehat{T(F)}$ and satisfy the property that
$$H(\pi_G(z\gamma),\psi)=\psi(z)H(\pi_G(\gamma),\psi)$$
for all $z\in Z_G(F)$, $\gamma\in G(F)^{\text{\emph{rss}}}$ and $\psi\in\widehat{T(F)}$ and let $\widetilde{H}\in C^\infty(\mathfrak{A}_{\widetilde{G}}\times \widehat{\widetilde{T}(F)})$ be defined via $\widetilde{H}(za_G,\widetilde{\psi})=\widetilde{\psi}(z)H(a_G,\widetilde{\psi}|_{T(F)})$ and $\widetilde{H}(a_{\widetilde{G}},\widetilde{\psi})=0$ if $a_{\widetilde{G}}\notin\mathfrak{A}_{\widetilde{G}}^+$. Then we have that $\mathcal{P}^{\widetilde{H}}(za_G)=\delta_z\otimes\mathcal{P}^H(a_G)$ and $\mathcal{P}^{\widetilde{H}}(a_{\widetilde{G}})=0$ if $a_{\widetilde{G}}\notin\mathfrak{A}_{\widetilde{G}}^+$. Moreover, we have $\widetilde{\lambda^H}=\lambda^{\widetilde{H}}\in\mathcal{D}(\mathfrak{A}_{\widetilde{G}}\times \widetilde{T},\Delta_{\det})$ with $\mathfrak{L}_{\widetilde{\lambda^H}}(\widetilde{\psi})=\widetilde{H}(\bullet,\widetilde{\psi})d\eta_{\widetilde{G}}$ as well as that, for $f=f^Z\otimes f^+\in C_c^\infty(\mathfrak{A}_{\widetilde{G}}^+)$ and $g\in C_c^\infty(\widetilde{T}(F))$ we have
$$\left<\widetilde{\lambda^H},f\otimes g\right>=\int_{Z_{\widetilde{G}}(F)}f^Z(z)\left<\mathfrak{R}_{\lambda^H}(f^+),\mathcal{M}_z(g)\right>d\mu_{Z_{\widetilde{G}}}$$
and
$$\mathfrak{R}_{\widetilde{\lambda^H}}(f^Z\otimes f^+)=(f^Z\otimes \mathfrak{R}_{\lambda^H}(f^+))^{Z_{\widetilde{G}}(F)T(F)}d\mu_{\widetilde{T}}.$$
\end{lem}
\begin{proof}
For $za_G\in \mathfrak{A}_{\widetilde{G}}^+$ and $g\in C_c^\infty(\widetilde{T}(F))$ we have by Fourier Inversion on $\widetilde{T}(F)/T(F)$ that
\begin{align*}
\left<\mathcal{P}^{\widetilde{H}}(za_G),g\right>=&\int_{\widehat{\widetilde{T}(F)}}H(a_G,\widetilde{\psi}|_{T(F)})\int_{\widetilde{T}(F)}g(y)\widetilde{\psi}(zy^{-1})d\mu_{\widetilde{T}}d\nu_{\widetilde{T}}\\
=&\int_{\widehat{T(F)}}H(a_G,\psi)\int_{T(F)^\perp}\int_{\widetilde{T}(F)/T(F)}\widetilde{\psi}(z\dot{y}^{-1})\int_{T(F)}g(\dot{y}t)\psi(t^{-1})d\mu_T\frac{d\mu_{\widetilde{T}}}{d\mu_T}d\nu_{\widetilde{T}}d\nu_T\\
=&\int_{\widehat{T(F)}}H(a_G,\psi)\int_{T(F)}g(zt)\psi(t^{-1})d\mu_Td\nu_T\\
=&\left<\mathcal{P}^H(a_G),(\mathcal{M}_zg)|_{T(F)}\right>\\
=&\left<\delta_z\otimes \mathcal{P}^H(a_G),g\right>
\end{align*}
and where we have $\mathcal{P}^{\widetilde{H}}(a_{\widetilde{G}})=0$ for $a_{\widetilde{G}}\notin\mathfrak{A}_{\widetilde{G}}^+$ by the definition of $\widetilde{H}$; this establishes the first claim. 

For $f=f^Z\otimes f^+\in C_c^\infty(\mathfrak{A}_{\widetilde{G}})$ supported on $\mathfrak{A}_{\widetilde{G}}^+$ we have by Fubini's Theorem that
\begin{align*}
\left<\widetilde{\lambda^H},f\otimes g\right>=&\int_{Z_{\widetilde{G}}(F)}f^Z(z)\int_{\mathfrak{A}_G}f^+(a_G)\int_{\widehat{T(F)}}H(a_G,\psi)\int_{T(F)}g(zt)\psi(t^{-1})d\mu_Td\nu_Td\eta_Gd\mu_{Z_{\widetilde{G}}}\\
=&\int_{\mathfrak{A}_G}f^+(a_G)\int_{\widehat{T(F)}}H(a_G,\psi)\int_{T(F)}\int_{Z_{\widetilde{G}}(F)}f^Z(z)g(zt)d\mu_{Z_{\widetilde{G}}}\psi(t^{-1})d\mu_Td\nu_Td\eta_G
\end{align*}
which shows that $\lambda^{\widetilde{H}}\in\mathcal{D}(\mathfrak{A}_{\widetilde{G}}\times \widetilde{T},\Delta_{\det})$ since the inner integral $\int_{Z_{\widetilde{G}}(F)}f^Z(z)g(zt)d\mu_{Z_{\widetilde{G}}}$ will vanish if $f\otimes g$ is supported away from $\Delta_{\det}$. Moreover, for $\widetilde{\psi}_0\in \widehat{\widetilde{T}(F)}$, $\psi_0=\widetilde{\psi}_0|_{T(F)}$ and $\widetilde{\psi}_{0,m}=\widetilde{\psi}_0\mathbf{1}_{\{y:|\text{ord}(\det(y))|\le m\}}$ we have by the calculation above that
\begin{align*}
&\left<\mathfrak{L}_{\widetilde{\lambda^H}}(\widetilde{\psi}),f^Z\otimes f^+\right>\\
=&\lim_{m\to\infty}\int_{Z_{\widetilde{G}}(F)}f^Z(z)\int_{\mathfrak{A}_G}f^+(a_G)\int_{\widehat{T(F)}}H(a_G,\psi)\int_{T(F)}\widetilde{\psi}_{0,m}(zt)\psi(t^{-1})d\mu_Td\nu_Td\eta_Gd\mu_{Z_{\widetilde{G}}}\\
=&\lim_{m\to\infty}\int_{Z_{\widetilde{G}}(F)}f^Z(z)\int_{\mathfrak{A}_G}f^+(a_G)\widetilde{\psi}_{0,m}(z)H(a_G,\psi_0)d\eta_Gd\mu_{Z_{\widetilde{G}}}\\
=&\int_{Z_{\widetilde{G}}(F)}f^Z(z)\int_{\mathfrak{A}_G}f^+(a_G)\widetilde{\psi}_0(z)H(a_G,\psi_0)d\eta_Gd\mu_{Z_{\widetilde{G}}}\\
=&\int_{\mathfrak{A}_{\widetilde{G}}}(f^Z\otimes f^+)(a_{\widetilde{G}})\widetilde{H}(a_{\widetilde{G}},\widetilde{\psi}_0)d\eta_{\widetilde{G}}
\end{align*}
as desired. The last statement regarding $\mathfrak{R}_{\lambda^{\widetilde{H}}}$ follows from the formula above since
\begin{align*}
\left<\widetilde{\lambda^H},f\otimes g\right>=&\int_{\mathfrak{A}_G}f^+(a_G)\int_{\widehat{T(F)}}H(a_G,\psi)\int_{T(F)}\int_{Z_{\widetilde{G}}(F)}f^Z(z)g(zt)d\mu_{Z_{\widetilde{G}}}\psi(t^{-1})d\mu_Td\nu_Td\eta_G\\
=&\int_{Z_{\widetilde{G}}(F)}f^Z(z)\int_{T(F)}\int_{\mathfrak{A}_G}f^+(a_G)\int_{\widehat{T(F)}}H(a_G,\psi)\psi(t^{-1})d\nu_Td\eta_Gg(zt)d\mu_Td\mu_{Z_{\widetilde{G}}}\\
=&\int_{Z_{\widetilde{G}}(F)}f^Z(z)\int_{T(F)}\mathfrak{L}_{\lambda^H}(f^+)(t)g(zt)d\mu_Td\mu_{Z_{\widetilde{G}}}\\
=&\int_{\widetilde{T}(F)}(f^Z\otimes \mathfrak{R}_{\lambda^H}(f^+))^{Z_{\widetilde{G}}(F)T(F)}(y)g(y)d\mu_{\widetilde{T}}
\end{align*}
and we are done.
\end{proof}

Let $\mathfrak{A}_{\widetilde{G}}^{\det(\widetilde{T})}\subset\mathfrak{A}_{\widetilde{G}}$ be the subset of elements $a_{\widetilde{G}}$ satisfying $\det(a_{\widetilde{G}})\in\det(\widetilde{T}(F))$. We define the distribution $\mathcal{V}(a_{\widetilde{G}})\in\mathcal{D}(\widetilde{T})$ via $\mathcal{V}(a_{\widetilde{G}})=0$ if $a_{\widetilde{G}}\notin \mathfrak{A}_{\widetilde{G}}^{\det(\widetilde{T})}$ and 
$$\left<\mathcal{V}(a_{\widetilde{G}}),g\right>=\int_{T(F)}g(y_{a_{\widetilde{G}}}t)d\mu_T$$
for $a_{\widetilde{G}}\in \mathfrak{A}_{\widetilde{G}}^{\det(\widetilde{T})}$ and any $y_{a_{\widetilde{G}}}\in\widetilde{T}(F)$ satisfying $\det(y_{a_{\widetilde{G}}})=\det(a_{\widetilde{G}})$. Similarly, for $\rho\in T(F)^\perp$ we define $\rho(a_{\widetilde{G}})=\rho(y_{a_{\widetilde{G}}})$ for any $y_{a_{\widetilde{G}}}$ satisfying $\det(y_{a_{\widetilde{G}}})=\det(a_{\widetilde{G}})$ and observe that this definition is independent of the particular choice of $y_{a_{\widetilde{G}}}$. With these notions in hand, we have the following.

\begin{lem}\label{ThingforA}
Let $h\in C^\infty(\mathfrak{A}_{\widetilde{G}}(F))$ be supported on $\mathfrak{A}_{\widetilde{G}}^{\det(\widetilde{T})}$ and $N\in C^\infty(\mathfrak{A}_{\widetilde{G}}\times \widetilde{T})$ be given by $N(a_{\widetilde{G}},\widetilde{\psi})=h(a_{\widetilde{G}})\mathbf{1}_{T(F)^\perp}(\widetilde{\psi})\widetilde{\psi}(a_{\widetilde{G}})$. Then $\mathcal{P}^N(a_{\widetilde{G}})=h(a_{\widetilde{G}})\mathcal{V}(a_{\widetilde{G}})$, for $\widetilde{\lambda}^{h,T(F)^\perp}=\lambda^N$ we have $\widetilde{\lambda}^{h,T(F)^\perp}\in \mathcal{D}(\mathfrak{A}_{\widetilde{G}}\times \widetilde{T},\Delta_{\det})$, $\mathfrak{L}_{\widetilde{\lambda}^{h,T(F)^\perp}}(\widetilde{\psi})=N(\bullet,\widetilde{\psi})d\eta_{\widetilde{G}}$ and for any $f\in C_c^\infty(\mathfrak{A}_{\widetilde{G}})$ we have
$$\mathfrak{R}_{\widetilde{\lambda}^{h,T(F)^\perp}}(f)=\left(\int_{\mathfrak{A}_G^{\det(\bullet)}}f(a_G^{\det(\bullet)})h(a_G^{\det(y)})d\eta_G^{\det(\bullet)}\right)d\mu_{\widetilde{T}}$$
where $y\mapsto \int_{\mathfrak{A}_G^{\det(y)}}f(a_G^{\det(y)})h(a_G^{\det(y)})d\eta_G^{\det(y)}$ is smooth and compactly supported.
\end{lem}
\begin{proof}
Letting $g\in C_c^\infty(\widetilde{T}(F))$ we have by Fourier Inversion on $\widetilde{T}(F)/T(F)$ that for $a_{\widetilde{G}}\in\mathfrak{A}_{\widetilde{G}}^{\det(\widetilde{T})}$ and $y_{a_{\widetilde{G}}}\in \widetilde{T}(F)$ with $\det(y_{a_{\widetilde{G}}})=\det(a_{\widetilde{G}})$ that
\begin{align*}
\left<\mathcal{P}^N(a_{\widetilde{G}}),g\right>=&h(a_{\widetilde{G}})\int_{T(F)^\perp}\mathcal{F}_{\widetilde{T}}(g)(\widetilde{\psi})\widetilde{\psi}(y_{a_{\widetilde{G}}})d\nu_{\widetilde{T}}\\
=&h(a_{\widetilde{G}})\int_{T(F)^\perp}\int_{\widetilde{T}(F)/T(F)}\widetilde{\psi}(y_{a_{\widetilde{G}}}\dot{t}^{-1})\int_{T(F)}g(\dot{t}s)d\mu_T\frac{d\mu_{\widetilde{T}}}{d\mu_T}d\nu_{\widetilde{T}}\\
=&h(a_{\widetilde{G}})\int_{T(F)}g(y_{a_{\widetilde{G}}}t)d\mu_T\\
=&\left<h(a_{\widetilde{G}})\mathcal{V}(a_{\widetilde{G}}),g\right>
\end{align*}
which establishes the first statement. We moreover see that
\begin{align*}
\left<\mathfrak{R}_{\lambda^{h,T(F)^\perp}}(f),g\right>=&\int_{\mathfrak{A}_{\widetilde{G}}^{\det(\widetilde{T})}}f(a_{\widetilde{G}})h(a_{\widetilde{G}})\int_{T(F)}g(y_{a_{\widetilde{G}}}t)d\mu_Td\eta_{\widetilde{G}}\\
=&\int_{\det(\widetilde{T}(F))}\int_{\mathfrak{A}_G^{\det(y)}}f(a_G^{\det(y)})h(a_G^{\det(y)})d\eta_G^{\det(y)}\int_{T(F)}g(\dot{y}t)d\mu_Td\mu_{F^\times}\\
=&\int_{\widetilde{T}(F)}\left(\int_{\mathfrak{A}_G^{\det(y)}}f(a_G^{\det(y)})h(a_G^{\det(y)})d\eta_G^{\det(y)}\right)g(y)d\mu_{\widetilde{T}}
\end{align*}
and the proof of the other statements is similar to the proofs of the analogous statements in Lemma \ref{GLdistdecomp}.

\end{proof}

We will also require the following which follows from a straightforward calculation.
\begin{lem}\label{GLdelta}
The map $\mathcal{P}_\delta:\mathfrak{A}_{\widetilde{G}}\to\mathcal{D}(\widetilde{T}(F))$ given by $\mathcal{P}_\delta(a_{\widetilde{G}})=\frac{1}{|D_{\widetilde{G}}(a_{\widetilde{G}})|^\frac{1}{2}}\sum_{t^\prime\in\widetilde{T}(F):\pi_{\widetilde{G}}(t^\prime)=a_{\widetilde{G}}}\delta_{t^\prime}$ is a smooth family with $\widetilde{\lambda}^\delta=\lambda^{\mathcal{P}_\delta}\in \mathcal{D}(\mathfrak{A}_{\widetilde{G}}\times \widetilde{T},\Delta_{\det})$ and where $\mathfrak{L}_{\widetilde{\lambda}^\delta}(\widetilde{\psi})=\frac{1}{|D_{\widetilde{G}}(a_{\widetilde{G}})|^\frac{1}{2}}\sum_{t^\prime\in\widetilde{T}(F):\pi_{\widetilde{G}}(t^\prime)=a_{\widetilde{G}}}\widetilde{\psi}(t^\prime)$ for all $\widetilde{\psi}\in\widehat{\widetilde{T}(F)}$ and $\mathfrak{R}_{\widetilde{\lambda}^\delta}(f)=\pi_{\widetilde{T}}^\ast(f)d\mu_{\widetilde{T}}$.
\end{lem}

\subsection{Computation of $\Theta_{\widetilde{\phi}}$}\label{GLfamilycomps}

We now use the various results established in $\S$\ref{helloGL} to compute $\Theta_{\widetilde{\phi}}$. The following is the computation of $\Theta_{\widetilde{\phi}}$ as well as the transfer operator $\mathfrak{S}_{\phi}^{\text{rss}}$ for functions with support on the regular semisimple locus; this answers in the affirmative Langlands' Question A for $\widetilde{\phi}:{^L\widetilde{T}}\to{^L\widetilde{G}}$. We remark that most of the work to prove Proposition \ref{GLQA} that remains simply involves properly organizing our results; most of the difficult work has already been done throughout $\S$\ref{helloGL}.
\begin{prop}\label{GLQA}
Let $\Theta_{\widetilde{\phi}}=\lambda^{H_{\widetilde{\phi}}}$ for $H_{\widetilde{\phi}}(a_{\widetilde{G}},\widetilde{\psi})=\chi_{\widetilde{\psi}}(a_{\widetilde{G}})$. Then we have that $\Theta_{\widetilde{\phi}}=\Theta_{\widetilde{\phi}}^++\Theta_{\widetilde{\phi}}^A+\Theta_{\widetilde{\phi}}^T$ for $\Theta_{\widetilde{\phi}}^+=\widetilde{\Theta_\phi^+}$, $\Theta_{\widetilde{\phi}}^A=\widetilde{\lambda}^{u_1,T(F)^\perp}$, and $\Theta_{\widetilde{\phi}}^T=\widetilde{\lambda}^\delta$. Furthermore, we have that $\Theta_{\widetilde{\phi}}\in\mathcal{D}(\mathfrak{A}_{\widetilde{G}}\times \widetilde{T},\Delta_{\det})$ with $\mathfrak{L}_{\Theta_{\widetilde{\phi}}}(\widetilde{\psi})=\chi_{\widetilde{\psi}}d\eta_{\widetilde{G}}$ for all $\widetilde{\psi}\in\widehat{\widetilde{T}(F)}$. Moreover, there exists an operator $\mathfrak{S}_{\widetilde{\phi}}^{\text{\emph{rss}}}:C_c^\infty(\mathfrak{A}_{\widetilde{G}})\to C_c^\infty(\widetilde{T}(F))$ satisfying $\mathfrak{R}_{\Theta_{\widetilde{\phi}}}(f)=\mathfrak{S}_{\widetilde{\phi}}^{\text{\emph{rss}}}(f)d\mu_{\widetilde{T}}$ for all $f\in C_c^\infty(\mathfrak{A}_{\widetilde{G}})$ where we have $\mathfrak{S}_{\widetilde{\phi}}^{\text{\emph{rss}}}=\widetilde{\mathfrak{S}}^++\widetilde{\mathfrak{S}}^A+\widetilde{\mathfrak{S}}^T$ with $\widetilde{\mathfrak{S}}^+$ supported on $\mathfrak{A}_{\widetilde{G}}^+$ with
$$\widetilde{\mathfrak{S}}^+(f^Z\otimes f^+)=(f^Z\otimes\mathfrak{S}^+(f^+))^{Z_{\widetilde{G}}(F)T(F)}$$
as well as
$$\widetilde{\mathfrak{S}}^A(f)(y)=\left(\int_{\mathfrak{A}_G^{\det(y)}}f(a_G^{\det(y)})u_1(a_G^{\det(y)})d\eta_G^{\det(y)}\right)$$
and
$$\widetilde{\mathfrak{S}}^T(f)=\pi_{\widetilde{T}}^\ast(f).$$
\end{prop}
\begin{proof}
By Proposition \ref{CharacterTable} we have that we may decompose
\begin{align*}
H_{\widetilde{\phi}}(a_{\widetilde{G}},\widetilde{\psi})=&H_{\widetilde{\phi}}(a_{\widetilde{G}},\widetilde{\psi})\mathbf{1}_{\mathfrak{A}_G^+\setminus \mathfrak{A}_{\widetilde{G}}^T}(a_{\widetilde{G}})\mathbf{1}_{\widehat{\widetilde{T}(F)}\setminus T(F)^\perp}(\widetilde{\psi})\\
&+H_{\widetilde{\phi}}(a_{\widetilde{G}},\widetilde{\psi})\mathbf{1}_{\mathfrak{A}_{\widetilde{G}}^T}(a_{\widetilde{G}})\\
&+H_{\widetilde{\phi}}(a_{\widetilde{G}},\widetilde{\psi})\mathbf{1}_{T(F)^\perp}(\widetilde{\psi}).
\end{align*}

Writing $H^{+\setminus T}(a_{\widetilde{G}},\widetilde{\psi})=H_{\widetilde{\phi}}(a_{\widetilde{G}},\widetilde{\psi})\mathbf{1}_{\mathfrak{A}_G^+\setminus \mathfrak{A}_{\widetilde{G}}^T}(a_{\widetilde{G}})\mathbf{1}_{\widehat{\widetilde{T}(F)}\setminus T(F)^\perp}(\widetilde{\psi})$ and $\lambda^{+\setminus T}=\lambda^{H^{+\setminus T}}$ we have by Lemma \ref{GLdistdecomp} that $\lambda^{+\setminus T}\in \mathcal{D}(\mathfrak{A}_{\widetilde{G}}\times \widetilde{T},\Delta_{\det})$ with $\mathfrak{L}_{\lambda^{+\setminus T}}(\widetilde{\psi})=\chi_{\widetilde{\psi}}\mathbf{1}_{\mathfrak{A}_G^+\setminus \mathfrak{A}_{\widetilde{G}}^{\widetilde{T}}}\mathbf{1}_{\widehat{\widetilde{T}(F)}\setminus T(F)^\perp}(\widetilde{\psi})d\eta_{\widetilde{G}}$ and that by Lemma \ref{GLdistdecomp} and Proposition \ref{SLtransferA} we have
$$\mathfrak{R}_{\lambda^{+\setminus T}}(f)=\widetilde{\mathfrak{S}}^+(f)d\mu_{\widetilde{T}}.$$

Writing $H^T(a_{\widetilde{G}},\widetilde{\psi})=H_{\widetilde{\phi}}(a_{\widetilde{G}},\widetilde{\psi})\mathbf{1}_{\mathfrak{A}_{\widetilde{G}}^T}(a_{\widetilde{G}})$ we similarly have by Lemmas \ref{GLdistdecomp} and \ref{GLdelta} that $\lambda^{H^T}\in \mathcal{D}(\mathfrak{A}_{\widetilde{G}}\times \widetilde{T},\Delta_{\det})$ with $\mathfrak{L}_{\lambda^{H^T}}(\widetilde{\psi})=\chi_{\widetilde{\psi}}\mathbf{1}_{\mathfrak{A}_{\widetilde{G}}^T}d\eta_{\widetilde{G}}$ and moreover that 
$$\mathfrak{R}_{\lambda^{H^T}}(f)=\left(\widetilde{\mathfrak{S}}^+(f)+\widetilde{\mathfrak{S}}^T(f)\right)d\mu_{\widetilde{T}}.$$

Writing $H^A(a_{\widetilde{G}},\widetilde{\psi})=H_{\widetilde{\phi}}(a_{\widetilde{G}},\widetilde{\psi})\mathbf{1}_{T(F)^\perp}(\widetilde{\psi})$ we have by Lemma \ref{ThingforA} that $\lambda^{H^A}\in \mathcal{D}(\mathfrak{A}_{\widetilde{G}}\times \widetilde{T},\Delta_{\det})$ with $\mathfrak{L}_{\lambda^{H^A}}(\widetilde{\psi})=\chi_{\widetilde{\psi}}\mathbf{1}_{T(F)^\perp}(\widetilde{\psi})d\eta_{\widetilde{G}}$ and
$$\mathfrak{R}_{\lambda^{H^A}}(f)=\widetilde{\mathfrak{S}}^A(f)d\mu_{\widetilde{T}}.$$
The result follows.
\end{proof}

\subsection{Stable Transfer for $\widetilde{G}(F)$}\label{GLtransfersec}

Similarly to the $G(F)$ case, for $g\in C^\infty(\widetilde{T}(F)^{\text{rss}})$ we define $\iota(g):T(F)\to \mathbb{C}$ to be given by $\iota(g)(t)=g(t)$ for $t\in \widetilde{T}(F)^{\text{rss}}$ and $\iota(g)(t)=0$ otherwise where furthermore, we set
$$C_{\mu_{\widetilde{T}}}^\infty(\widetilde{T}(F))=\left\{g\in C^\infty(\widetilde{T}(F)^{\text{rss}}):\text{supp}(\iota(g))\text{ is compactly supported, }\int_{\widetilde{T}(F)}|\iota(g)(t)|d\mu_{\widetilde{T}}<\infty\right\}$$
where, attaining the first inclusion by restricting functions to $\widetilde{T}(F)^{\text{rss}}$, we identify $C_c^\infty(\widetilde{T}(F))\subset C_{\mu_{\widetilde{T}}}^\infty(\widetilde{T}(F))\subset C^\infty(\widetilde{T}(F))$. We may now answer in the affirmative Langlands' Question B of \cite{SetT} for $\widetilde{\phi}$.
\begin{thm}\label{transferforGL}
There exists an operator $\mathfrak{S}_{\widetilde{\phi}}:C_b^\infty(\mathfrak{A}_{\widetilde{G}})\to C_{\mu_{\widetilde{T}}}^\infty(\widetilde{T})$ which satisfies 
$$\int_{\mathfrak{A}_{\widetilde{G}}}f(a_{\widetilde{G}})\chi_{\widetilde{\psi}}(a_{\widetilde{G}})d\eta_{\widetilde{G}}=\int_{\widetilde{T}(F)}\widetilde{\mathfrak{S}}_b(f)(t)\widetilde{\psi}(t)d\eta_{\widetilde{T}}$$
for all $\widetilde{\psi}\in\widehat{\widetilde{T}(F)}$ and that $\mathfrak{S}_{\widetilde{\phi}}|_{C_c^\infty(\mathfrak{A}_{\widetilde{G}}(F)^{\text{\emph{rss}}})}=\mathfrak{S}_{\widetilde{\phi}}^{\text{\emph{rss}}}$.
\end{thm}
\begin{proof}
For any $f\in C_b^\infty(\mathfrak{A}_{\widetilde{G}})$ we may decompose $f=f\mathbf{1}_{\mathfrak{A}_{\widetilde{G}}^+}+f\mathbf{1}_{\mathfrak{A}_{\widetilde{G}}\setminus\mathfrak{A}_{\widetilde{G}}^+}$ where $f\mathbf{1}_{\mathfrak{A}_{\widetilde{G}}\setminus\mathfrak{A}_{\widetilde{G}}^+}$ may be handled similarly to what was done in the proof of Theorem \ref{SLMainThm} by using the explicit formulas of Proposition \ref{GLQA}. We henceforth assume $f=f\mathbf{1}_{\mathfrak{A}_{\widetilde{G}}^+}$ is supported on $\mathfrak{A}_{\widetilde{G}}^+$. By Proposition \ref{GLproductOI} it further suffices to consider the case $f=f^Z\otimes f^+$ for $f^Z\in C_c^\infty(Z_{\widetilde{G}}(F))$ and $f^+\in C_b^\infty(\mathfrak{A}_G)$ supported on $\mathfrak{A}_G^+$. We have by Proposition \ref{GLQA} and Theorem \ref{SLMainThm} that the operator
$$\widetilde{\mathfrak{S}}_{\widetilde{\phi}}(f^Z\otimes f^+)=(f^Z\otimes\mathfrak{S}_\phi(f))^{Z_{\widetilde{G}}(F)T(F)}$$
satisfies the desired properties.
\end{proof}



\begin{thebibliography}{99}




\bibitem{ADeBMKT}
J. Adler, S. DeBacker, \emph{Murnaghan-Kirillov theory for supercuspidal representations of general linear groups}, J. Reine Angew. Math \textbf{575} (2004), pp. 1-35.

\bibitem{ADeBBT}
\bysame, \emph{Some applications of Bruhat-Tits theory to harmonic analysis on the Lie algebra of a reductive p-adic group}, with appendices by R. Huntsinger and G. Prasad, Michigan Math. J. \textbf{50} (2002), no, 2, pp. 263-286.



\bibitem{ASformula}
J. Adler, L. Spice, \emph{Supercuspidal characters of reductive {$p$}-adic groups}, Amer. J. Math. \textbf{131} (2009), no. 4, pp. 1137-1210.

\bibitem{SL2}
J. Adler, S. DeBacker, P. J. Sally, Jr., L. Spice, \emph{Supercuspidal characters of $\text{\emph{SL}}_2$ over a $p$-adic field}, in \emph{Harmonic analysis on reductive, p-adic groups}, R. Doran, P. Sally, L. Spice, eds., Contemporary Mathematics, vol. 543, American Mathematical Society, Providence, RI, 2011., pp. 19-70.




\bibitem{BT1}
F. Bruhat, J. Tits, \emph{Groupes r\'eductifs sur un corps local}, Inst. Hautes \'Etudes Sci. Publ. Math. (1972), no. 41, pp. 5-251.

\bibitem{BT2}
\bysame, \emph{Groupes r\'eductifs sur un corps local. {II}. {S}ch\'emas en groupes. {E}xistence d'une donn\'ee radicielle valu\'ee}, Inst. Hautes \'Etudes Sci. Publ. Math. (1984), no. 60, pp. 197-376.




\bibitem{Nilpotent}
D. Collingwood, W. McGovern, \emph{Nilpotent orbits in semisimple lie algebras}, Chapman and Hall, London, 1993.


\bibitem{DeThesis}
S. DeBacker, \emph{On supercuspidal characters of $\text{\emph{GL}}_\ell$, $\ell$ a prime}, ProQuest LLC, Ann Arbor, MI, 1997, Thesis (Ph.D.)-The University of Chicago.

\bibitem{DeBLCE}
\bysame, \emph{Homogeneity Results for Invariant Distributions of a Reductive $p$-adic Group}, Annales scientifiques de l'\'{E}cole Normale Sup\'{e}rieure, Serie 4, vol. 35 (2002) no. 3, pp. 391-422.


\bibitem{DeBS}
S. DeBacker, L. Spice, \emph{Stability of character sums for positive-depth supercuspidal representations}, arXiv:1310.3306.



\bibitem{FLN}
E. Frenkel, R. Langlands, B.C. Ngo, \emph{Formule des Trace et Fonctorialit\'{e}: le D\'{e}but d'un Programme}, Ann. Math. Qu\'{e}. \textbf{34}, (2010).

\bibitem{SLLpack}
S. Gelbart, A. Knapp, \emph{L-Indistinguishability and R Groups for the Special Linear Group}, Adv. in Math. \textbf{43}, pp. 101-121 (1982).

\bibitem{GGPS}
I. M. Gelfand, M. I. Graev, I. I. Piatetski-Shapiro, \emph{Representation theory and automorphic functions}, Saunders, 1968.

  
\bibitem{MurnHowe}
J. Hakim, F. Murnaghan, \emph{Distinguished tame supercuspidal representations}, Int. Math. Res. Pap. IMRP 166(2) (2008).



\bibitem{Howe}
R. Howe, \emph{Tamely ramified supercuspidal representations of {$\text{GL}_{n}$}}, Pacific J. Math. \textbf{73} (1977), no. 2, pp. 437-460.

\bibitem{MyThesis}
D. Johnstone, \emph{A Gelfand-Graev formula and stable transfer factors for $\text{SL}_n(F)$}, ProQuest LLC, Ann Arbor, MI, 2017, Thesis (Ph.D.)-The University of Chicago.

\bibitem{ThesisExtension}
\bysame, \emph{Stable Transfer for Maximal Elliptic Tori in $\text{SL}_n(F)$ and $\text{GL}_n(F)$}, in preparation.

\bibitem{Symn}
D. Johnstone, Z. Luo, \emph{On the Stable Transfer for $\text{Sym}^{\text{n}}$ Lifting of $\text{GL}_2$}, arXiv:2002.09551.

\bibitem{Symn2}
\bysame, in preparation.



\bibitem{Exhaustion}
J-L. Kim, \emph{Supercuspidal representations: an exhaustion theorem}, J. Amer. Math. Soc. \textbf{20} (2007), no. 2, pp. 273-320.
  
\bibitem{KimMurn}
J-L. Kim, F. Murnaghan, \emph{Character expansions and unrefined minimal K-types}, Amer. J. Math. \textbf{125} (2003), pp. 1199-1234.
  


\bibitem{KNotes}
R. Kottwitz, \emph{Harmonic Analysis on Reductive p-adic Groups and Lie Algebras}, in \emph{Harmonic Analysis, the Trace Formula, and Shimura Varieties} (2003), pp. 393-522.



\bibitem{BE}
R. P. Langlands, \emph{Beyond Endoscopy}, http://publications.ias.edu.

\bibitem{SetT}
\bysame, \emph{Singularit\'{e}s et transfert}, Ann. Math. Qu\'{e}. \textbf{37} (2013), no. 2, pp. 173-253.

\bibitem{LS}
R. P. Langlands, D. Shelstad, \emph{On the Definition of Transfer Factors}, Math. Ann. \textbf{278} (1987), pp. 219-271.


\bibitem{MP1}
A. Moy, G. Prasad, \emph{Unrefined minimal {$K$}-types for {$p$}-adic groups}, Invent. Math. \textbf{116} (1994), no. 1-3, pp. 393-408.

\bibitem{MP2}
\bysame, \emph{Jacquet functors and unrefined minimal {$K$}-types}, Comment. Math. Helv. \textbf{71} (1996), no. 1, pp. 98-121.

\bibitem{LCE}
F. Murnaghan, \emph{Local character expansions and Shalika germs for GL(n)}, Math. Ann. \textbf{304} (1996), pp. 423-455.



\bibitem{Rao}
R. Rao, \emph{Orbital Integrals in Reductive Groups}, Annals of Math. \textbf{96} (1972) no. 3, pp. 505-510.

\bibitem{Rudin}
W. Rudin, \emph{Fourier Analysis on Groups}, Hoboken, NJ: Wiley; 2011.

\bibitem{YiannisGG}
Y. Sakellaridis, \emph{Transfer operators and Hankel transforms between relative trace formulas, II: Rankin--Selberg theory},  arXiv:1805.04640.

\bibitem{SS} 
P. J. Sally, Jr., J. A. Shalika, \emph{Characters of the discrete series of representations of $\text{\emph{SL}}(2)$ over a local field}, Proc. Nat. Acad. Sci. U.S.A. \textbf{61} (1968), 1231-1237.

\bibitem{Scholze}
P. Scholze, \emph{The local Langlands correspondence for $\text{GL}_n$ over p-adic fields}, Invent. Math. \textbf{192} (2013), no. 3, pp. 663-715.

\bibitem{Distr}
L. Schwartz, \emph{Th\'{e}orie des distributions}, \textbf{1-2}, Hermann (1951).

\bibitem{Shalika}
J. A. Shalika, \emph{A Theorem on Semi-Simple P-adic Groups}, Annals of Math. \textbf{95} (1972) no. 3, pp. 226-242.

\bibitem{SLL} L. Spice, \emph{Supercuspidal characters of $\text{\emph{SL}}_\ell$ over a $p$-adic field, $\ell$ a prime}, Amer. J. Math. \textbf{127} (2005) pp. 51-100.


\bibitem{SpiceNew} \bysame, \emph{Explicit asymptotic expansions for tame supercuspidal characters}, arXiv:1701.02417.


\bibitem{Takahashi}
T. Takahashi, \emph{Characters of cuspidal unramified series for central simple algebras of prime degree}, J. Math. Kyoto Univ. \textbf{32-4} (1992), pp. 873-888.

\bibitem{Titscor}
J. Tits, \emph{Reductive groups over local fields}, in \emph{Automorphic forms, representations and L-functions (Proc. Sympos. Pure Math., Oregon State Univ., Corvallis, Ore., 1977), Part 1}, Providence, R.I., 1979, pp. 29-69.

\bibitem{vanDijk}
G. van Dijk, \emph{Computation of certain induced characters of $\mathfrak{p}$-adic groups}, Math. Ann. \textbf{199} (1972), pp. 229-240.

\bibitem{Vigneras}
M-F Vign\'{e}ras, \emph{Caract\'{e}risation des int\'{e}grales orbitales sur un groupe r\'{e}ductif $p$-adique}, Journal of the Faculty of Science, University of Tokyo, Sect IA, Vol 28, $n^o$ 3, pp. 945-961 (1982). 

\bibitem{JKYBT}
J-K. Yu, \emph{Bruhat-Tits theory and buildings}, in \emph{Ottawa Lectures on Admissible Representations of Reductive p-adic Groups}, pp. 53-79.

\bibitem{LocLang}
\bysame, \emph{On the Local Langlands Correspondence for Tori}, in \emph{Ottawa Lectures on Admissible Representations of Reductive p-adic Groups}, pp. 177-183.

\bibitem{Construction}
\bysame, \emph{Construction of tame supercuspidal representations}, J. Amer. Math. Soc. \textbf{14} (2001), no. 3, pp. 579-622.



\end{thebibliography}
\end{document}